\documentclass[10pt]{amsart}
\usepackage{amscd}
\usepackage{amssymb}
\usepackage[all]{xy}
\usepackage{graphicx}
\author{Amnon Yekutieli}
\title{Mixed Resolutions and Simplicial Sections} 
\date{8 December 2005}
\address{Department of  Mathematics,
Ben Gurion University, Be'er Sheva 84105, Israel}
\email{amyekut@math.bgu.ac.il}
\thanks{
Primary: 14F10; Secondary: 18G30, 16E45, 18G10.}
\thanks{ This work was partially supported by the 
US - Israel Binational Science Foundation}

\newtheorem{thm}[equation]{Theorem}
\newtheorem{cor}[equation]{Corollary}
\newtheorem{prop}[equation]{Proposition}
\newtheorem{lem}[equation]{Lemma}
\theoremstyle{definition}
\newtheorem{dfn}[equation]{Definition}
\newtheorem{rem}[equation]{Remark}
\newtheorem{exa}[equation]{Example}

\numberwithin{equation}{section}

\newcommand{\iso}{\xrightarrow{\simeq}}
\newcommand{\inj}{\hookrightarrow}

\newcommand{\xar}{\xrightarrow}
\newcommand{\opn}{\operatorname}
\newcommand{\cat}[1]{\operatorname{\mathsf{#1}}}

\newcommand{\rmitem}[1]{\item[\text{\textup{(#1)}}]}
\newcommand{\mfrak}[1]{\mathfrak{#1}}
\newcommand{\mcal}[1]{\mathcal{#1}}
\newcommand{\msf}[1]{\mathsf{#1}}
\newcommand{\mbf}[1]{\mathbf{#1}}
\newcommand{\mrm}[1]{\mathrm{#1}}
\newcommand{\mbb}[1]{\mathbb{#1}}

\newcommand{\tup}[1]{\textup{#1}}
\newcommand{\bsym}[1]{\boldsymbol{#1}}
\newcommand{\boplus}{\bigoplus\nolimits}

\newcommand{\til}[1]{\tilde{#1}}
\newcommand{\what}[1]{\widehat{#1}}

\newcommand{\set}[1]{\{ #1 \}}

\newcommand{\K}{\mbb{K} \hspace{0.05em}}
\renewcommand{\d}{\mrm{d}}

\newcommand{\bprod}{{\textstyle \prod}}

\newcommand{\hatotimes}[1]{\, \what{\otimes}_{#1} \,}

\begin{document}

\begin{abstract}
We introduce the notions of mixed resolutions and simplicial 
sections, and prove a theorem relating them. This result is used 
(in another paper) 
to study deformation quantization in algebraic geometry. 
\end{abstract}

\maketitle

\setcounter{section}{-1}
\section{Introduction}

Let $\K$ be a field of characteristic $0$. In this paper we 
present several technical results about the geometry of 
$\K$-schemes. These results were discovered in 
the course of work on deformation quantization in algebraic 
geometry, and they play a crucial role in \cite{Ye3}. This role 
will be explained at the end of the introduction. The idea behind 
the constructions in this paper can be traced back to old work of 
Bott \cite{Bo, HY}. 

Let $\pi : Z \to X$ be a morphism of $\K$-schemes, and let 
$\bsym{U} = \{ U_{(0)}, \ldots, U_{(m)} \}$ be an open covering of 
$X$. A {\em simplicial section} $\bsym{\sigma}$ of $\pi $, based on 
the covering $\bsym{U}$, consists of a family of morphisms 
$\sigma_{\bsym{i}} : \bsym{\Delta}^q_{\K} \times 
U_{\bsym{i}} \to Z$, 
where $\bsym{i} = (i_0, \ldots, i_q)$ is a multi-index; 
$\bsym{\Delta}^q_{\K}$ is the $q$-dimensional geometric simplex; 
and 
$U_{\bsym{i}} := U_{(i_{0})} \cap \cdots \cap U_{(i_{q})}$. 
The morphisms $\sigma_{\bsym{i}}$ are required to be compatible 
with $\pi$ and to satisfy simplicial relations. 
See Definition \ref{dfn3.2} for details.
An important example of a simplicial section is mentioned at the 
end of the introduction. 

Another notion we introduce is that of {\em mixed resolution}. 
Here we assume the $\K$-scheme $X$ is smooth and separated,
and each of the open sets $U_{(i)}$ in the covering 
$\bsym{U}$ is affine. Given a quasi-coherent 
$\mcal{O}_X$-module $\mcal{M}$ we define its mixed resolution
$\opn{Mix}_{\bsym{U}}^{}(\mcal{M})$. This is a complex of sheaves 
on $X$, concentrated in non-negative degrees. 
As the name suggests, this 
resolution mixes two distinct types of resolutions: a de Rham type 
resolution which is related to the sheaf $\mcal{P}_X$ of principal 
parts of $X$ and its Grothendieck connection, and a 
simplicial-\v{C}ech 
type resolution which is related to the covering
$\bsym{U}$. The precise definition is too complicated to state 
here -- see Section \ref{sec4}. 

Let $\msf{C}^{+}(\cat{QCoh} \mcal{O}_X)$ denote the abelian 
category of bounded below complexes of quasi-coherent 
$\mcal{O}_X$-modules. For any 
$\mcal{M} \in \msf{C}^{+}(\cat{QCoh} \mcal{O}_X)$ the mixed 
resolution 
$\opn{Mix}_{\bsym{U}}^{}(\mcal{M})$ is defined by totalizing the 
double complex 
$\boplus_{p,q} \opn{Mix}^q_{\bsym{U}}(\mcal{M}^p)$.
The derived category of $\K$-modules is denoted by 
$\msf{D}(\cat{Mod} \K)$.

\begin{thm} \label{0.3}
Let $X$ be a smooth separated $\K$-scheme, and let 
$\bsym{U} = \linebreak \{ U_{(0)}, \ldots, U_{(m)} \}$ be an affine
open covering of $X$.
\begin{enumerate}
\item There is a functorial quasi-isomorphism
$\mcal{M} \to \opn{Mix}^{}_{\bsym{U}}(\mcal{M})$
for $\mcal{M} \in \linebreak \msf{C}^{+}(\cat{QCoh} \mcal{O}_X)$.
\item Given 
$\mcal{M} \in \msf{C}^{+}(\cat{QCoh} \mcal{O}_X)$,
the canonical morphism 
$\Gamma \big( X, \opn{Mix}^{}_{\bsym{U}}(\mcal{M}) \big) \to
\mrm{R} \Gamma \bigl( X, \opn{Mix}^{}_{\bsym{U}}(\mcal{M}) \big)$
in $\msf{D}(\cat{Mod} \K)$ is an isomorphism.
\item The quasi-isomorphism in part \tup{(1)} induces a 
functorial isomorphism \linebreak
$\Gamma \big( X, \opn{Mix}^{}_{\bsym{U}}(\mcal{M}) \big) \cong
\mrm{R} \Gamma(X, \mcal{M})$ 
in $\msf{D}(\cat{Mod} \K)$.
\end{enumerate}
\end{thm}

This is repeated as Theorem \ref{thm4.1} in the body of the paper. 
Note that part (3) is a formal consequence of parts (1) and (2).

A useful corollary of the theorem is the following (see Corollary 
\ref{cor4.1}). Suppose $\mcal{M}$ and $\mcal{N}$ are two 
complexes in $\msf{C}^{+}(\cat{QCoh} \mcal{O}_X)$, and 
$\phi : \opn{Mix}^{}_{\bsym{U}}(\mcal{M}) \to 
\opn{Mix}^{}_{\bsym{U}}(\mcal{N})$ 
is a $\K$-linear quasi-isomorphism. Then 
\[ \Gamma(X, \phi) : \Gamma \bigl( X, \opn{Mix}^{}_{\bsym{U}} 
\mcal{M}) \bigr) 
\to \Gamma \bigl( X, \opn{Mix}^{}_{\bsym{U}}(\mcal{N}) \bigr) \]
is a quasi-isomorphism. 

Here is the connection between simplicial sections and 
mixed resolutions. 

\begin{thm} \label{thm0.4}
Let $X$ be a smooth separated $\K$-scheme, let
$\pi : Z \to X$ be a morphism of schemes, and let 
$\bsym{U}$ be an affine open covering of $X$. 
Suppose $\bsym{\sigma}$ is a simplicial section of $\pi$ 
based on $\bsym{U}$. Let 
$\mcal{M}_1, \ldots, \mcal{M}_r, \mcal{N}$ be quasi-coherent 
$\mcal{O}_{X}$-modules, and let
\[ \phi : \prod_{i = 1}^r  \pi^{\what{*}}
(\mcal{P}_X \otimes_{\mcal{O}_{X}} \mcal{M}_i) 
\to \pi^{\what{*}}
(\mcal{P}_X \otimes_{\mcal{O}_{X}} \mcal{N}) \]
be a continuous $\mcal{O}_Z$-multilinear sheaf morphism on $Z$. 
Then there is an induced 
$\K$-multilinear sheaf morphism 
\[ \bsym{\sigma}^*(\phi) :
\prod_{i = 1}^r \opn{Mix}_{\bsym{U}}(\mcal{M}_i) \to 
\opn{Mix}_{\bsym{U}}(\mcal{N})  \]
on $X$.
\end{thm}

In the theorem, the continuity and the complete pullback 
$\pi^{\what{*}}$ 
refer to the dir-inv structures on these sheaves, which are 
explained in Section \ref{sec1}. 
A more detailed statement is Theorem \ref{thm3.3}
in the body of the paper. 

Let us explain, in vague terms, how Theorem \ref{thm0.4}, or rather
Theorem \ref{thm3.3}, is used in the paper \cite{Ye3}. Let $X$ be a 
smooth separated $n$-dimensional $\K$-scheme. As we know from 
the work of Kontsevich \cite{Ko}, there are two important sheaves 
of DG Lie algebras on $X$, namely the sheaf 
$\mcal{T}^{}_{\mrm{poly}, X}$ of poly derivations, and the sheaf 
$\mcal{D}^{}_{\mrm{poly}, X}$ of poly differential operators.
Suppose $\bsym{U}$ is some affine open covering of $X$.
The inclusions
$\mcal{T}^{}_{\mrm{poly}, X} \to
\opn{Mix}^{}_{\bsym{U}}(\mcal{T}^{}_{\mrm{poly}, X})$
and 
$\mcal{D}^{}_{\mrm{poly}, X} \to
\opn{Mix}^{}_{\bsym{U}}(\mcal{D}^{}_{\mrm{poly}, X})$
are then quasi-isomorphisms of sheaves of DG Lie algebras 
(cf.\ Theorem \ref{0.3}). The goal is to find an $\mrm{L}_{\infty}$ 
quasi-isomorphism 
\[ \Psi : \opn{Mix}^{}_{\bsym{U}}(\mcal{T}^{}_{\mrm{poly}, X})
\to \opn{Mix}^{}_{\bsym{U}}(\mcal{D}^{}_{\mrm{poly}, X}) \]
between these sheaves of DG Lie algebras. Having such 
an $\mrm{L}_{\infty}$ quasi-isomorphism pretty much
implies the solution of the deformation quantization problem for 
$X$. 

Let $\opn{Coor} X$ denote the coordinate bundle of $X$. This is an 
infinite dimensional bundle over $X$,
endowed with an action of the group $\mrm{GL}_n(\K)$. 
Let $\opn{LCC} X$ be the quotient bundle
$\opn{Coor} X\, /\, \mrm{GL}_n(\K)$. 
In \cite{Ye4} we proved that if the covering $\bsym{U}$ is fine 
enough (the condition is that each open set $U_{(i)}$ admits an 
\'etale morphism to $\mbf{A}^n_{\K}$), then the projection 
$\pi : \opn{LCC} X \to X$ admits a simplicial section 
$\bsym{\sigma}$. 

Now the universal deformation formula of Kontsevich \cite{Ko} 
gives rise to a continuous $\mrm{L}_{\infty}$ quasi-isomorphism 
\[ \mcal{U} : \pi^{\what{*}}
(\mcal{P}_X \otimes_{\mcal{O}_{X}} \mcal{T}^{}_{\mrm{poly}, X})
\to \pi^{\what{*}}
(\mcal{P}_X \otimes_{\mcal{O}_{X}} \mcal{D}^{}_{\mrm{poly}, X}) \]
on $\opn{LCC} X$. This means that there is a sequence of 
continuous $\mcal{O}_{\opn{LCC} X}$-multilinear sheaf morphisms 
\[ \mcal{U}_r : \bprod^r\,  \pi^{\what{*}}
(\mcal{P}_X \otimes_{\mcal{O}_{X}} \mcal{T}^{}_{\mrm{poly}, X})
\to \pi^{\what{*}}
(\mcal{P}_X \otimes_{\mcal{O}_{X}} \mcal{D}^{}_{\mrm{poly}, X}) , \]
$r \geq 1$, satisfying very complicated identities. 
Using Theorem \ref{thm3.3} 
we obtain a sequence of multilinear sheaf morphisms 
\[ \bsym{\sigma}^*(\mcal{U}_r) :
\bprod^r \opn{Mix}^{}_{\bsym{U}}(\mcal{T}^{}_{\mrm{poly}, X})
\to \opn{Mix}^{}_{\bsym{U}}(\mcal{D}^{}_{\mrm{poly}, X})  \]
on $X$. 
After twisting these morphisms suitably (this is needed due to the 
presence of the Grothendieck connection; cf.\ \cite{Ye2}) 
we obtain the desired $\mrm{L}_{\infty}$ quasi-isomorphism $\Psi$.

We believe that mixed resolutions, and the results of this paper, 
shall have additional applications in algebraic geometry (e.g.\ 
algebro-geometric versions of results on index theorems in 
differential geometry, cf.\ \cite{NT}; or a proof of Kontsevich's 
famous yet unproved claim on Hochschild cohomology of a scheme
\cite[Claim 8.4]{Ko}).

\section{Review of Dir-Inv Modules}
\label{sec1}

We begin the paper with a review of the concept of dir-inv 
structure, which was introduced in \cite{Ye2}. 
A dir-inv structure is a generalization of adic topology. 

Let $C$ be a commutative ring.
We denote by $\cat{Mod} C$ the category of $C$-modules. 

\begin{dfn}
\begin{enumerate}
\item Let $M \in \cat{Mod} C$. An {\em inv module 
structure} on $M$ is an inverse system 
$\{ \mrm{F}^i M \}_{i \in \mbb{N}}$ of $C$-submodules
of $M$. The pair 
$(M, \{ \mrm{F}^i M \}_{i \in \mbb{N}})$
is called an {\em inv $C$-module}.
\item Let $(M, \{ \mrm{F}^i M \}_{i \in \mbb{N}})$ and
$(N, \{ \mrm{F}^i N \}_{i \in \mbb{N}})$ be two inv 
$C$-modules.
A function $\phi : M \to N$ 
($C$-linear or not) is said to be {\em continuous}
if for every $i \in \mbb{N}$ there exists $i' \in \mbb{N}$ such 
that $\phi(\mrm{F}^{i'} M) \subset \mrm{F}^i N$.
\item Define $\cat{Inv} \cat{Mod} C$
to be the category whose objects are the inv $C$-modules, 
and whose morphisms are the continuous $C$-linear 
homomorphisms.  
\end{enumerate}
\end{dfn}

There is a full and faithful embedding of categories
$\cat{Mod} C \inj \cat{Inv} \cat{Mod} C$,
$M \mapsto (M, \{ \ldots, 0, 0 \})$.

Recall that a directed set is a partially ordered set $J$ 
with the property that for any $j_1, j_2 \in J$ there exists 
$j_3 \in J$ such that $j_1, j_2 \leq j_3$.

\begin{dfn}
\begin{enumerate}
\item Let $M \in \cat{Mod} C$. A {\em 
dir-inv module structure} on $M$ is a direct system 
$\{ \mrm{F}_j M \}_{j \in J}$ of $C$-submodules 
of $M$, indexed by a nonempty directed set $J$, 
together with an inv module structure on 
each $\mrm{F}_j M$, such that for every $j_1 \leq j_2$ the 
inclusion $\mrm{F}_{j_1} M \inj \mrm{F}_{j_2} M$
is continuous. The pair 
$(M, \{ \mrm{F}_j M \}_{j \in J})$
is called a {\em dir-inv $C$-module}.
\item Let $(M, \{ \mrm{F}_j M \})_{j \in J}$ and
$(N, \{ \mrm{F}_k N \}_{k \in K})$ 
be two dir-inv $C$-modules.
A function $\phi : M \to N$ 
($C$-linear or not)  is said to be {\em continuous}
if for every $j \in J$ there exists $k \in K$ 
such that $\phi(\mrm{F}_j M) \subset \mrm{F}_k N$, and
$\phi : \mrm{F}_j M \to \mrm{F}_k N$ is a continuous 
homomorphism between these two inv $C$-modules.
\item Define $\cat{Dir} \cat{Inv} \cat{Mod} C$
to be the category whose objects are the dir-inv $C$-modules, 
and whose morphisms are the continuous $C$-linear 
homomorphisms.  
\end{enumerate}
\end{dfn}

An inv $C$-module $M$ can be endowed with a
dir-inv module structure $\{ \mrm{F}_j M \}_{j \in J}$, where 
$J := \{ 0 \}$ and $\mrm{F}_0 M := M$.
Thus we get a full and faithful embedding 
$\cat{Inv} \cat{Mod} C \inj 
\cat{Dir} \cat{Inv} \cat{Mod} C$. 

Inv modules and dir-inv modules come in a few ``flavors'': 
trivial, discrete and complete. 
A {\em discrete inv module} is one which is isomorphic, in
$\cat{Inv} \cat{Mod} C$, to an object of 
$\cat{Mod} C$ (via the canonical embedding above). 
A {\em complete inv module} is an inv module 
$(M, \{ \mrm{F}^i M \}_{i \in \mbb{N}})$
such that the canonical map 
$M \to \lim_{\leftarrow i} \mrm{F}^i M$
is bijective. A {\em discrete} (resp.\ {\em complete})
{\em dir-inv module} 
is one which is isomorphic, in 
$\cat{Dir} \cat{Inv} \cat{Mod} C$, to a dir-inv module
$(M, \{ \mrm{F}_j M \}_{j \in J})$, where all the inv 
modules $\mrm{F}_j M$ are discrete (resp.\ complete), 
and the canonical map
$\lim_{j \to} \mrm{F}_j M \to M$
in $\cat{Mod} C$ is bijective. A {\em trivial dir-inv 
module} is one which is isomorphic to an object of 
$\cat{Mod} C$. Discrete dir-inv modules are complete, but 
there are also other complete modules, as the next example shows.

\begin{exa} \label{exa2.4}
Assume $C$ is noetherian and $\mfrak{c}$-adically complete 
for some ideal $\mfrak{c}$. Let $M$ be a finitely 
generated $C$-module, and define 
$\mrm{F}^i M := \mfrak{c}^{i+1} M$. Then 
$\{ \mrm{F}^i M \}_{i \in \mbb{N}}$ is called the 
{\em $\mfrak{c}$-adic inv structure}, and of course 
$(M, \{ \mrm{F}^i M \}_{i \in \mbb{N}})$ is a complete inv module.
Next consider an arbitrary $C$-module $M$. We take
$\{ \mrm{F}_j M \}_{j \in J}$ to be the collection of 
finitely generated $C$-submodules of $M$. 
This dir-inv module structure on $M$ 
is called the {\em $\mfrak{c}$-adic dir-inv structure}. 
Again $(M, \{ \mrm{F}_j M \}_{j \in J})$ is a complete 
dir-inv $C$-module. Note that a finitely generated
$C$-module $M$ is discrete as inv module iff 
$\mfrak{c}^{i} M = 0$ for $i \gg 0$; and a $C$-module
is discrete as dir-inv module iff it is a direct 
limit of discrete finitely generated modules.
\end{exa}

The category $\cat{Dir} \cat{Inv} \cat{Mod} C$ is additive.
Given a collection $\{ M_k \}_{k \in K}$ of dir-inv modules, the 
direct sum $\boplus_{k \in K} M_k$ has a structure of dir-inv 
module, making it into the coproduct of 
$\{ M_k \}_{k \in K}$ in the category
$\cat{Dir} \cat{Inv} \cat{Mod} C$. Note that if the index set
$K$ is infinite and each $M_k$ is a nonzero 
discrete inv module, then $\boplus_{k \in K} M_k$ is a 
discrete dir-inv module which is not trivial.
The tensor product 
$M \otimes_{C} N$
of two dir-inv modules is again a dir-inv module.
There is a completion functor
$M \mapsto \what{M}$. 
(Warning: if $M$ is complete then 
$\what{M} = M$, but
it is not known if $\what{M}$ is complete for arbitrary 
$M$.) The completed tensor product is  
$M \what{\otimes}_{C} N :=
\what{M \otimes_{C} N}$. 
Completion commutes with direct sums: if 
$M \cong \boplus_{k \in K} M_k$ then 
$\what{M} \cong  \boplus_{k \in K} \what{M}_k$.
See \cite{Ye2} for full details. 

A {\em graded dir-inv module} (or graded object in 
$\cat{Dir} \cat{Inv} \cat{Mod} C$)
is a direct sum $M = \boplus_{k \in \mbb{Z}} M_k$, 
where each $M_k$ is a dir-inv module. A {\em DG algebra}
in $\cat{Dir} \cat{Inv} \cat{Mod} C$ is a graded dir-inv 
module $A = \boplus_{k \in \mbb{Z}} A^k$, together 
with continuous $C$-(bi)linear functions 
$\mu : A \times A \to A$ and $\d : A \to A$, 
which make $A$ into a DG $C$-algebra.
If $A$ is a super-commutative associative unital DG algebra in 
$\cat{Dir} \cat{Inv} \cat{Mod} C$, and $\mfrak{g}$ is a DG Lie 
Algebra in $\cat{Dir} \cat{Inv} \cat{Mod} C$, then 
$A \hatotimes{C} \mfrak{g}$ is a DG Lie 
Algebra in $\cat{Dir} \cat{Inv} \cat{Mod} C$. 

Let $A$ be a super-commutative associative unital DG algebra in 
$\cat{Dir} \cat{Inv} \cat{Mod} C$. A {\em DG $A$-module} 
in $\cat{Dir} \cat{Inv} \cat{Mod} C$ is a graded object
$M$ in $\cat{Dir} \cat{Inv} \cat{Mod} C$,
together with continuous $C$-(bi)linear functions 
$\mu : A \times M \to M$ and $\d : M \to M$, 
which make $M$ into a DG $A$-module in the usual sense. 
A {\em DG $A$-module Lie algebra} in
$\cat{Dir} \cat{Inv} \cat{Mod} C$ is a DG Lie algebra
$\mfrak{g}$ in $\cat{Dir} \cat{Inv} \cat{Mod} C$,
together with a continuous $C$-bilinear function
$\mu : A \times \mfrak{g} \to \mfrak{g}$, such that
such that $\mfrak{g}$ becomes a DG $A$-module, and 
\[ [a_1 \gamma_1, a_2 \gamma_2] = (-1)^{i_2 j_1}
a_1 a_2 \,[\gamma_1, \gamma_2] \]
for all $a_k \in A^{i_k}$ and 
$\gamma_k \in \mfrak{g}^{j_k}$.

All the constructions above can be geometrized. 
Let $(Y, \mcal{O})$ be a commutative ringed space
over $\K$, i.e.\ $Y$ is a topological space, and $\mcal{O}$ is 
a sheaf of commutative $\K$-algebras on $Y$. 
We denote by $\cat{Mod} \mcal{O}$ the category of 
$\mcal{O}$-modules on $Y$. 

\begin{exa} \label{exa2.1}
Geometrizing Example \ref{exa2.4}, let $\mfrak{X}$ be a noetherian 
formal scheme, with defining ideal $\mcal{I}$. Then any coherent 
$\mcal{O}_{\mfrak{X}}$-module $\mcal{M}$ is an 
inv $\mcal{O}_{\mfrak{X}}$-module, with system of submodules
$\{ \mcal{I}^{i+1} \mcal{M} \}_{i \in \mbb{N}}$,
and $\mcal{M} \cong \what{\mcal{M}}$; cf.\ 
\cite{EGA-I}. We call an $\mcal{O}_{\mfrak{X}}$-module
{\em dir-coherent} if it is the direct limit of coherent 
$\mcal{O}_{\mfrak{X}}$-modules. Any dir-coherent module 
is quasi-coherent, but it is not known if the converse is true.
At any rate, a dir-coherent $\mcal{O}_{\mfrak{X}}$-module 
$\mcal{M}$ is a dir-inv $\mcal{O}_{\mfrak{X}}$-module, where we 
take $\{ \mrm{F}_j \mcal{M} \}_{j \in J}$ to be the collection of 
coherent submodules of $\mcal{M}$. Any 
dir-coherent $\mcal{O}_{\mfrak{X}}$-module is then a complete 
dir-inv module. This dir-inv module structure on $\mcal{M}$ 
is called the {\em $\mcal{I}$-adic dir-inv structure}. Note that 
a coherent $\mcal{O}_{\mfrak{X}}$-module $\mcal{M}$ 
is discrete as inv module iff 
$\mcal{I}^{i} \mcal{M} = 0$ for $i \gg 0$; and a 
dir-coherent $\mcal{O}_{\mfrak{X}}$-module is discrete as dir-inv 
module iff it is a direct limit of discrete coherent modules.
\end{exa}

If $f : (Y', \mcal{O}') \to (Y, \mcal{O})$ is a morphism of 
ringed spaces and
$\mcal{M} \in \cat{Dir} \cat{Inv} \cat{Mod} \mcal{O}$,
then there is an obvious structure of dir-inv 
$\mcal{O}'$ -module on $f^{*} \mcal{M}$, and we define
$f^{\what{*}} \mcal{M} := \what{f^* \mcal{M}}$. 
If $\mcal{M}$ is a graded object in
$\cat{Dir} \cat{Inv} \cat{Mod} \mcal{O}$,
then the inverse images $f^* \mcal{M}$
and $f^{\what{*}} \mcal{M}$ are graded objects in 
$\cat{Dir} \cat{Inv} \cat{Mod} \mcal{O}'$.
If $\mcal{G}$ is an algebra (resp.\ a DG algebra) in 
$\cat{Dir} \cat{Inv} \cat{Mod} \mcal{O}$, then 
$f^* \mcal{G}$ and  $f^{\what{*}}\, \mcal{G}$ are
algebras (resp.\ DG algebras) in 
$\cat{Dir} \cat{Inv} \cat{Mod} \mcal{O}'$. 
Given $\mcal{N} \in \cat{Dir} \cat{Inv} \cat{Mod} \mcal{O}'$
there is an obvious 
dir-inv $\mcal{O}$-module structure on $f_* \mcal{N}$.

\begin{exa} \label{exa2.3}
Let $(Y, \mcal{O})$ be a ringed space and $V \subset Y$ an open 
set. For a dir-inv $\mcal{O}$-module $\mcal{M}$ there is an 
obvious way to make $\Gamma(V, \mcal{M})$ into a dir-inv 
$\mcal{O}$-module. If $\mcal{M}$ is a complete inv $\mcal{O}$-module
then $\Gamma(V, \mcal{M})$ is a complete inv $\mcal{O}$-module.
If $V$ is quasi-compact and $\mcal{M}$ 
is a complete dir-inv $\mcal{O}$-module, then 
$\Gamma(V, \mcal{M})$ is a complete dir-inv
$\Gamma(V, \mcal{O})$-module.
\end{exa}

\section{Complete Thom-Sullivan Cochains}
\label{sec2}

 From here on $\K$ is a field of characteristic $0$.
Let us begin with some abstract notions about cosimplicial modules 
and their normalizations, following \cite{HS} and \cite{HY}. 
We use the notation $\cat{Mod} \K$ and $\cat{DGMod} \K$ for the 
categories of $\K$-modules and DG (differential graded) 
$\K$-modules respectively. 

Let $\bsym{\Delta}$ denote the category with objects the ordered 
sets $[q] := \{ 0, 1, \ldots, q \}$, $q \in \mbb{N}$. The 
morphisms $[p] \to [q]$ are the order preserving functions,
and we write
$\bsym{\Delta}^q_p := \opn{Hom}_{\bsym{\Delta}}([p], [q])$. 
The $i$-th co-face map 
$\partial^i : [p] \to [p + 1]$
is the injective function that does not take the value $i$;
and the $i$-th co-degeneracy map
$\mrm{s}^i : [p] \to [p - 1]$
is the surjective function that takes the value $i$ twice.
All morphisms in $\bsym{\Delta}$ are compositions of various 
$\partial^i$ and $\mrm{s}^i$. 

An element of $\bsym{\Delta}^q_p$ may be thought of as a sequence
$\bsym{i} = (i_0, \ldots, i_p)$ of integers with
$0 \leq i_0 \leq \cdots \leq i_p \leq q$. 
Given $\bsym{i} \in \bsym{\Delta}^m_q$, 
$\bsym{j} \in \bsym{\Delta}_m^p$ and
$\alpha \in \bsym{\Delta}^q_p$, 
we sometimes write 
$\alpha_*(\bsym{i}) := \bsym{i} \circ \alpha \in \bsym{\Delta}^m_p$ 
and 
$\alpha^*(\bsym{j}) := \alpha  \circ\bsym{j} \in 
\bsym{\Delta}_m^q$.

Let $\cat{C}$ be some category. A {\em cosimplicial object} in 
$\cat{C}$ is a functor $C : \bsym{\Delta} \to \cat{C}$.
We shall usually refer to the cosimplicial object as
$C = \set{C^p}_{p \in \mbb{N}}$,
and for any $\alpha \in \bsym{\Delta}_p^q$ the corresponding 
morphism in $\cat{C}$ will be denoted by 
$\alpha^* : C^p \to C^q$.
A {\em simplicial object} in 
$\cat{C}$ is a functor $C : \bsym{\Delta}^{\mrm{op}} \to \cat{C}$.
The notation for a simplicial object
will be $C = \set{C_p}_{p \in \mbb{N}}$ and
$\alpha_* : C_q \to C_p$.

Suppose $M = \set{M^q}_{q \in \mbb{N}}$ is a cosimplicial 
$\K$-module. The {\em standard normalization} of $M$
is the DG module $\mrm{N} M$ defined as follows:
$\mrm{N}^q M := \bigcap_{i = 0}^{q - 1} \opn{Ker}
(\mrm{s}^i : M^q \to M^{q - 1})$.
The differential is
$\partial := \sum_{i = 0}^{q+1} (-1)^i \partial^i :
\mrm{N}^q M \to \mrm{N}^{q + 1} M$.
We get a functor
$\mrm{N} : \bsym{\Delta} \cat{Mod} \K \to \cat{DGMod} \K$.

For any $q$ let $\bsym{\Delta}^q_{\K}$ be the {\em geometric 
$q$-dimensional simplex}
\[ \bsym{\Delta}^q_{\K} := \mrm{Spec}\, 
\K [t_0, \ldots, t_q] / (t_0 + \cdots + t_q - 1) . \] 
The $i$-th vertex of $\bsym{\Delta}^q_{\K}$ is the $\K$-rational 
point $x$ such that $t_i(x) = 1$ and $t_j(x) = 0$ for all 
$j \neq i$. 
We identify the vertices of $\bsym{\Delta}^q_{\K}$ with the 
ordered set $[q] = \{ 0, 1, \ldots, q \}$. 
For any $\alpha : [p] \to [q]$ in $\bsym{\Delta}$ there is a 
unique linear morphism
$\alpha: \bsym{\Delta}^p_{\K} \to \bsym{\Delta}^q_{\K}$
extending it, and in this way 
$\{ \bsym{\Delta}^q_{\K} \}_{q \in \mbb{N}}$ 
is a cosimplicial scheme.

For a $\K$-scheme $X$ we write 
$\Omega^{p}(X) :=$ 
$\Gamma(X, \Omega^{p}_{X / \K})$. 
Taking $X := \bsym{\Delta}^q_{\K}$ we have a super-commutative
associative unital DG $\K$-algebra
$\Omega^{}(\bsym{\Delta}^q_{\K}) = 
\boplus_{p \in \mbb{N}} \Omega^{p}(\bsym{\Delta}^q_{\K})$, 
that is generated as $\K$-algebra by the elements
$t_0, \ldots, t_q, \mrm{d} t_0, \ldots, \mrm{d} t_q$. 
The collection
$\{ \Omega^{}(\bsym{\Delta}^q_{\K}) \}_{q \in \mbb{N}}$ 
is a simplicial DG algebra, namely a functor from 
$\bsym{\Delta}^{\mrm{op}}$ to the category of DG $\K$-algebras.

In \cite{HY} we made use of the Thom-Sullivan normalization 
$\til{\mrm{N}}^{} M$ of a cosimplicial $\K$-module $M$. For some 
applications (specifically \cite{Ye3}) a complete version of this 
construction is needed. Recall that for 
$M, N \in \cat{Dir} \cat{Inv} \cat{Mod} \K$
we can define the complete tensor product
$N \, \what{\otimes} \, M$.
The $\K$-modules $\Omega^{q}(\bsym{\Delta}^{l}_{\K})$
are always considered as discrete inv modules, so
$\Omega(\bsym{\Delta}^{l}_{\K})$ is a discrete dir-inv DG 
$\K$-algebra.

\begin{dfn} \label{dfn3.1}
Suppose $M= \{ M^{q} \}_{q \in \mbb{N}}$ is a cosimplicial 
dir-inv $\K$-module, namely each 
$M^q \in \cat{Dir} \cat{Inv} \cat{Mod} \K$, and 
the morphisms $\alpha^* : M^p \to M^q$, for 
$\alpha \in \bsym{\Delta}^q_p$, are continuous $\K$-linear 
homomorphisms. Let
\begin{equation} \label{eqn5.4}
\what{\til{\mrm{N}}}{}^{q} M \subset
\prod_{l=0}^{\infty} \left( 
\Omega^{q}(\bsym{\Delta}^{l}_{\K}) \, \what{\otimes} \, 
M^{l} \right)
\end{equation}
be the submodule consisting of all sequences
$(u_{0}, u_{1}, \ldots)$, 
with $u_l \in \Omega^{q}(\bsym{\Delta}^{l}_{\K}) 
\, \what{\otimes} \, M^{l}$,
such that
\begin{equation} \label{eqn4.1}
(\bsym{1} \otimes \alpha^*)(u_k) = 
(\alpha_* \otimes \bsym{1})(u_l) \in 
\Omega^{q}(\bsym{\Delta}^{k}_{\K}) \, \what{\otimes} \,  M^{l}
\end{equation}
for all $k, l \in \mbb{N}$ and all $\alpha \in \bsym{\Delta}_k^l$.
Define a coboundary operator
$\partial : \what{\til{\mrm{N}}}{}^{q} M \to 
\what{\til{\mrm{N}}}{}^{q + 1} M$
using the exterior derivative 
$\mrm{d} : \Omega^{q}(\bsym{\Delta}^{l}_{\K}) \to
\Omega^{q + 1}(\bsym{\Delta}^{l}_{\K})$.
The resulting DG $\K$-module
$\big( \what{\til{\mrm{N}}} M, \partial \big)$ is called the 
{\em complete Thom-Sullivan normalization of $M$}. 
\end{dfn}

The $\K$-module 
$\what{\til{\mrm{N}}} M = \boplus_{q \in \mbb{N}} 
\what{\til{\mrm{N}}}{}^q M$
is viewed as an abstract module. We obtain a functor
\[ \what{\til{\mrm{N}}} : 
\bsym{\Delta} \cat{Dir} \cat{Inv} \cat{Mod} \K 
\to \cat{DGMod} \K . \]

\begin{rem}
In case each $M^l$ is a discrete dir-inv module one has 
$\Omega^{q}(\bsym{\Delta}^{l}_{\K}) \, \what{\otimes} \,  M^{l}
= \Omega^{q}(\bsym{\Delta}^{l}_{\K}) \otimes M^{l}$,
and therefore  
$\what{\til{\mrm{N}}}{} M = \til{\mrm{N}} M$. 
\end{rem}

The standard normalization $\mrm{N} M$ also makes sense here, via 
the forgetful functor
$\bsym{\Delta} \cat{Dir} \cat{Inv} \cat{Mod} \K \to 
\bsym{\Delta} \cat{Mod} \K$. 
The two normalizations $\what{\til{\mrm{N}}}$ and
$\mrm{N}$ are related as follows. Let
$\int_{\bsym{\Delta}^l} : \Omega^{}(\bsym{\Delta}^l_{\K})
\to \K$
be the $\K$-linear map of degree $-l$
defined by integration on the compact real 
$l$-dimensional simplex, namely 
$\int_{\bsym{\Delta}^l} \mrm{d} t_1 \wedge \cdots \wedge
\mrm{d} t_l = \frac{1}{l !}$ etc. 
Suppose each dir-inv module $M^l$ is complete, 
so that using \cite[Proposition 1.5]{Ye2} 
we get a functorial $\K$-linear homomorphism
\[ \int_{\bsym{\Delta}^l} : 
\Omega^{}(\bsym{\Delta}^{l}_{\K}) \, \what{\otimes} \, 
M^{l} \to \K \, \what{\otimes} \, M^{l} \cong M^l . \]

\begin{prop} \label{prop4.4}
Suppose $M= \{ M^{q} \}_{q \in \mbb{N}}$ is a cosimplicial 
dir-inv $\K$-module, with all dir-inv modules $M^q$ complete. 
Then the homomorphisms $\int_{\bsym{\Delta}^l}$ induce 
a quasi-isomorphism 
\[ \int_{\bsym{\Delta}} : \what{\til{\mrm{N}}} M \to \mrm{N} M  \]
in $\cat{DGMod} \K$.
\end{prop}

\begin{proof}
This is a complete version of \cite[Theorem 1.12]{HY}.
Let $\bsym{\Delta}^{l}$ be the simplicial set 
$\bsym{\Delta}^{l} := \opn{Hom}_{\bsym{\Delta}}(-, [l])$;
so its set of $p$-simplices is $\bsym{\Delta}^{l}_p$. 
Define $C_l$ to 
be the algebra of normalized cochains on $\bsym{\Delta}^{l}$, 
namely
\[ C_l := \mrm{N} \opn{Hom}_{\cat{Sets}}(\bsym{\Delta}^{l}, \K)
\cong \opn{Hom}_{\cat{Sets}}(\bsym{\Delta}^{l, \mrm{nd}}, \K) . \]
Here $\bsym{\Delta}^{l, \mrm{nd}}$ is the (finite) set of 
nondegenerate simplices, i.e.\ those sequences
$\bsym{i} = (i_0, \ldots, i_p)$ satisfying
$0 \leq i_0 < \cdots < i_p \leq l$. 
As explained in \cite[Appendix A]{HY} we have simplicial DG 
algebras $C = \{ C_l \}_{l \in \mbb{N}}$ and 
$\Omega^{}(\bsym{\Delta}^{}_{\K}) =
\{ \Omega^{}(\bsym{\Delta}^{l}_{\K}) \}_{l \in \mbb{N}}$, 
and a homomorphism of simplicial DG modules
$\rho : \Omega^{}(\bsym{\Delta}^{}_{\K}) \to C$.

It turns out (this is work of Bousfield-Gugenheim)
that $\rho$ is a homotopy equivalence in 
$\bsym{\Delta}^{\mrm{op}} \cat{DGMod} \K$, 
i.e.\ there are simplicial homomorphisms 
$\phi : C \to \Omega^{}(\bsym{\Delta}^{}_{\K})$,
$h : C \to C$ and 
$h' : \Omega^{}(\bsym{\Delta}^{}_{\K}) \to
\Omega^{}(\bsym{\Delta}^{}_{\K})$ 
such that 
$\bsym{1} - \rho \circ \phi = h \circ \mrm{d} + \mrm{d} \circ h$
and
$\bsym{1} - \phi \circ \rho = h' \circ \mrm{d} + \mrm{d} \circ 
h'$.

Now for 
$M = \{ M^q \} \in 
\bsym{\Delta} \cat{Dir} \cat{Inv} \cat{Mod} \K$
and
$N = \{ N_q \} \in \bsym{\Delta}^{\mrm{op}} \cat{Mod} \K$
let $N \, \what{\otimes}_{\leftarrow} \, M$
be the complete version of \cite[formula (A.1)]{HY}, so that in 
particular
$\Omega^{}(\bsym{\Delta}^{}_{\K}) \, \what{\otimes}_{\leftarrow} \,
M \cong \what{\til{\mrm{N}}} M$
and
$C \, \what{\otimes}_{\leftarrow} \, M \cong \mrm{N} M$.
Moreover 
\[ \rho \, \what{\otimes}_{\leftarrow} \, \bsym{1}_M =
\int_{\bsym{\Delta}} : \what{\til{\mrm{N}}} M \to \mrm{N} M . \]
It follows that $\int_{\bsym{\Delta}}$ is a homotopy 
equivalence in $\cat{DGMod} \K$.
\end{proof}

Suppose $A = \{ A^q \}_{q \in \mbb{N}}$ is a cosimplicial DG 
algebra in $\cat{Dir} \cat{Inv} \cat{Mod} \K$
(not necessarily associative nor commutative).
This is a pretty complicated object: for every $q$ we have a DG 
algebra $A^q = \boplus_{i \in \mbb{Z}} \, A^{q, i}$ in
$\cat{Dir} \cat{Inv} \cat{Mod} \K$. For every
$\alpha \in \bsym{\Delta}^q_p$ there is a continuous
DG algebra homomorphism $\alpha^* : A^p \to A^q$,
and the $\alpha^*$ have to satisfy the 
simplicial relations.

Anyhow, both $\what{\til{\mrm{N}}} A$ and $\mrm{N} A$ are 
DG algebras. For $\what{\til{\mrm{N}}} A$ the DG algebra structure 
comes from that of the DG algebras
$\Omega(\bsym{\Delta}^{l}_{\K}) \, \what{\otimes} \, A^{l}$,
via the embeddings (\ref{eqn5.4}).
In case each $A^l$ is an associative super-commutative 
unital DG $\K$-algebra, then so is $\what{\til{\mrm{N}}} A$.
Likewise for DG Lie algebras. 
(The algebra $\mrm{N} A$, with its Alexander-Whitney product, 
is very noncommutative.)

Assume that each $A^{q, i}$ is complete, so that the 
integral 
$\int_{\bsym{\Delta}} : \what{\til{\mrm{N}}} A \to \mrm{N} A$ 
is defined. This 
is not a DG algebra homomorphism. However:

\begin{prop}
Suppose $A = \{ A^{q} \}_{q \in \mbb{N}}$ is a cosimplicial DG 
algebra in \linebreak
$\cat{Dir} \cat{Inv} \cat{Mod} \K$, with all 
$A^{q}$ complete. Then the homomorphisms 
$\int_{\bsym{\Delta}^l}$ induce an isomorphism of graded algebras
\[ \mrm{H}(\int_{\bsym{\Delta}}) : 
\mrm{H} \what{\til{\mrm{N}}} A \iso \mrm{H} \mrm{N} A . \]
\end{prop}

\begin{proof}
This is a complete variant of \cite[Theorem 1.13]{HY}. The proof 
is identical, after replacing ``$\otimes$'' with 
``$\what{\otimes}$'' where needed; cf.\ proof of previous 
proposition.
\end{proof}

\begin{rem}
If $A$ is associative then presumably $\int_{\bsym{\Delta}}$ 
extends to an $\mrm{A}_{\infty}$ quasi-isomorphism
$\what{\til{\mrm{N}}} A \to \mrm{N} A$.
\end{rem}

\section{Commutative \v{C}ech Resolutions}
\label{sec3}

In this section $\K$ is a field of characteristic $0$
and $X$ is a noetherian topological space. We denote by $\K_X$ the 
constant sheaf $\K$ on $X$. We will be interested in the category 
$\cat{Dir} \cat{Inv} \cat{Mod} \K_X$, whose objects are sheaves of 
$\K$-modules on $X$ with dir-inv structures. 
Note that any open set $V \subset X$ is quasi-compact. 

Let $X = \bigcup_{i = 0}^m U_{(i)}$ be an open covering, 
which we denote by $\bsym{U}$. For any 
$\bsym{i} = (i_0, \ldots, i_q) \in \bsym{\Delta}_q^m$ define
$U_{\bsym{i}} := U_{(i_0)} \cap \cdots \cap U_{(i_q)}$,
and let $g_{\bsym{i}} : U_{\bsym{i}} \to X$ be the inclusion.
Given a dir-inv $\K_X$-module $\mcal{M}$ and natural number $q$
we define a sheaf
\[ \mrm{C}^q(\bsym{U}, \mcal{M}) :=
\prod_{\bsym{i} \in \bsym{\Delta}_q^m} 
g_{\bsym{i} *}\, g_{\bsym{i}}^{-1}\, \mcal{M} . \]
This is a finite product. 
For an open set $V \subset X$ we then have
\[ \Gamma \bigl( V, \mrm{C}^q(\bsym{U}, \mcal{M}) \bigr) = 
\prod_{\bsym{i} \in \bsym{\Delta}_q^m}
\Gamma(V \cap U_{\bsym{i}}, \mcal{M}) . \] 
For any $\bsym{i}$ the $\K$-module 
$\Gamma(V \cap U_{\bsym{i}}, \mcal{M})$
has a dir-inv structure. Hence 
$\Gamma(V, \mrm{C}^q(\bsym{U}, \mcal{M}))$
is a dir-inv $\K$-module. If $\mcal{M}$ happens to be a
complete dir-inv $\K_X$-module then 
$\Gamma(V, \mrm{C}^q(\bsym{U}, \mcal{M}))$
is a complete dir-inv $\K$-module, since each 
$V \cap U_{\bsym{i}}$ is quasi-compact. 

Keeping $V$ fixed we get a cosimplicial dir-inv $\K$-module
$\big\{ \Gamma(V, \mrm{C}^q(\bsym{U}, \mcal{M})) 
\big\} _{q \in \mbb{N}}$.
Applying the functors $\mrm{N}^q$ and 
$\what{\til{\mrm{N}}}{}^q$ we obtain $\K$-modules
$\mrm{N}^q \Gamma(V, \mrm{C}(\bsym{U}, \mcal{M}))$
and
$\what{\til{\mrm{N}}}{}^q \Gamma(V, \mrm{C}(\bsym{U}, \mcal{M}))$.
As we vary $V$ these become pre\-sheaves of 
$\K$-modules, which we denote by
$\mrm{N}^q \mrm{C} (\bsym{U}, \mcal{M})$ and
$\what{\til{\mrm{N}}}{}^q \mrm{C} (\bsym{U}, \mcal{M})$.

Recall that a simplex $\bsym{i} = (i_0, \ldots, i_q)$ 
is nondegenerate if $i_0 < \cdots < i_q$.
Let $\bsym{\Delta}_q^{m, \mrm{nd}}$ be the set of non-degenerate 
simplices inside $\bsym{\Delta}_q^m$.

\begin{lem} \label{lem3.5}
For every $q$ the presheaves
\[ \mrm{N}^q \mrm{C} (\bsym{U}, \mcal{M}) : V \mapsto
\mrm{N}^q \Gamma \bigl( V, \mrm{C}(\bsym{U}, \mcal{M}) \bigr) \]
and
\[ \what{\til{\mrm{N}}}{}^q \mrm{C} (\bsym{U}, \mcal{M}) : V \mapsto
\what{\til{\mrm{N}}}{}^q \Gamma \bigl( 
V, \mrm{C}(\bsym{U}, \mcal{M}) \bigr) \]
are sheaves. There is a functorial isomorphism of sheaves
\begin{equation} \label{eqn3.2}
\mrm{N}^{q} \mrm{C}(\bsym{U}, \mcal{M})
\cong
{\displaystyle \prod_{\bsym{i} \in \bsym{\Delta}_q^{m, \mrm{nd}}}} 
\,
g_{\bsym{i} *}\, g_{\bsym{i}}^{-1}\, \mcal{M} ,
\end{equation}
and functorial embeddings of sheaves
\begin{equation} \label{eqn3.1}
\what{\til{\mrm{N}}}{}^{q} \mrm{C}(\bsym{U}, \mcal{M})
\inj 
{\displaystyle \prod_{l \in \mbb{N}}}\ \
{\displaystyle \prod_{\bsym{i} \in 
\bsym{\Delta}_l^{m}}} \,
g_{\bsym{i} *}\, g_{\bsym{i}}^{-1}\,
\bigl( \Omega^{q}(\bsym{\Delta}^l_{\K})
\, \what{\otimes} \, \mcal{M} \bigr) 
\end{equation}
and
\begin{equation} \label{eqn3.3}
\mcal{M} \inj 
\what{\til{\mrm{N}}}{}^{0} \mrm{C}(\bsym{U}, \mcal{M}) .
\end{equation} 
\end{lem}

\begin{proof}
Since $\{ \mrm{C}^q(\bsym{U}, \mcal{M}) \}_{q \in \mbb{N}}$
is a cosimplicial sheaf we get the isomorphism (\ref{eqn3.2}).

As for $\what{\til{\mrm{N}}}{}^{q} \mrm{C}(\bsym{U}, \mcal{M})$,
consider the sheaf 
$\Omega^{q}(\bsym{\Delta}^l_{\K}) \, \what{\otimes} \, \mcal{M}$
on $X$. Take any open set $V \subset X$ and 
$\bsym{i} \in \bsym{\Delta}_q^{m}$. Since $V \cap U_{\bsym{i}}$
is quasi-compact we have 
\[ \begin{aligned}
\Omega^{q}(\bsym{\Delta}^l_{\K}) \, \what{\otimes} \, 
\Gamma(V \cap U_{\bsym{i}}, \mcal{M}) & \cong
\Gamma \big( V \cap U_{\bsym{i}},
\Omega^{q}(\bsym{\Delta}^l_{\K}) \, \what{\otimes} \, 
\mcal{M} \big) \\
& = \Gamma \big( V, g_{\bsym{i} *}\, g_{\bsym{i}}^{-1}\,
(\Omega^{q}(\bsym{\Delta}^l_{\K}) \, \what{\otimes} \, 
\mcal{M}) \big) .
\end{aligned} \]
By Definition \ref{dfn3.1} there is an exact sequence of 
presheaves on $X$:
\[ \begin{aligned}
0 \to \what{\til{\mrm{N}}}{}^{q} \mrm{C}(\bsym{U}, \mcal{M})
& \to 
{\displaystyle \prod_{l \in \mbb{N}}}\ \
{\displaystyle \prod_{\bsym{i} \in 
\bsym{\Delta}_l^{m}}} \,
g_{\bsym{i} *}\, g_{\bsym{i}}^{-1}\,
\bigl( \Omega^{q}(\bsym{\Delta}^l_{\K})
\, \what{\otimes} \, \mcal{M} \bigr) \\
& \xar{\bsym{1} \otimes \alpha^* - \alpha_* \otimes \bsym{1}}
{\displaystyle \prod_{k, l \in \mbb{N}}}\ \
{\displaystyle \prod_{\alpha \in 
\bsym{\Delta}_k^{l}}} \
{\displaystyle \prod_{\bsym{i} \in 
\bsym{\Delta}_l^{m}}} \,
g_{\bsym{i} *}\, g_{\bsym{i}}^{-1}\,
\bigl( \Omega^{q}(\bsym{\Delta}^k_{\K})
\, \what{\otimes} \, \mcal{M} \bigr) .
\end{aligned} \]
Since the presheaves in the middle and on the right are actually 
sheaves, it follows that 
$\what{\til{\mrm{N}}}{}^{q} \mrm{C}(\bsym{U}, \mcal{M})$
is also a sheaf. 

Finally the embedding (\ref{eqn3.3}) comes from the embeddings
$\mcal{M} \inj \Omega^{0}(\bsym{\Delta}^l_{\K})
\, \what{\otimes} \, \mcal{M}$,
$w \mapsto 1 \otimes w$.
\end{proof}

Thus we have complexes of sheaves
$\mrm{N} \mrm{C}(\bsym{U}, \mcal{M})$ and
$\what{\til{\mrm{N}}} \mrm{C}(\bsym{U}, \mcal{M})$.
There are functorial homomorphisms
$\mcal{M} \to \mrm{N} \mrm{C}(\bsym{U}, \mcal{M})$
and
$\mcal{M} \to \what{\til{\mrm{N}}} \mrm{C}(\bsym{U}, \mcal{M})$.
Note that the complex 
$\Gamma(X, \mrm{N} \mrm{C}(\bsym{U}, \mcal{M}))$ 
is nothing but the usual global \v{C}ech complex of $\mcal{M}$ for 
the covering $\bsym{U}$.

\begin{dfn}
The complex 
$\what{\til{\mrm{N}}} \mrm{C}(\bsym{U}, \mcal{M})$
 is called the 
{\em commutative \v{C}ech resolution} of $\mcal{M}$.
\end{dfn}

The reason for the name is that 
$\what{\til{\mrm{N}}} \mrm{C}(\bsym{U}, \mcal{O}_X)$
is a sheaf of super-commutative DG algebras, as can be seen from 
the next lemma.

\begin{lem} \label{lem3.2}
Suppose $\mcal{M}_1, \ldots, \mcal{M}_r, \mcal{N}$ 
are dir-inv $\K_X$-modules, and 
$q_1, \ldots, q_r \in \mbb{N}$. Let $q := q_1 + \cdots + q_r$.
Suppose that for every 
$l \in \mbb{N}$ and 
$\bsym{i} \in \bsym{\Delta}_l^{m}$
we are given $\K$-multilinear sheaf maps
\[ \begin{aligned}
\phi_{q_1, \ldots, q_r, \bsym{i}} & : 
\big( \Omega^{q_1}(\bsym{\Delta}^{l}_{\K}) \, \what{\otimes} \,
(\mcal{M}_1|_{U_{\bsym{i}}}) \big) \times \cdots \times
\big( \Omega^{q_r}(\bsym{\Delta}^{l}_{\K}) \, \what{\otimes} \,
(\mcal{M}_r|_{U_{\bsym{i}}}) \big) \\
& \qquad \longrightarrow
\Omega^{q}(\bsym{\Delta}^{l}_{\K}) \, \what{\otimes} \,
(\mcal{N}|_{U_{\bsym{i}}}) 
\end{aligned} \]
that are continuous \tup{(}for the dir-inv module structures\tup{)}, 
and are compatible with the simplicial 
structure as in Definition \tup{\ref{dfn3.1}}.
Then there are unique $\K$-multilinear sheaf maps
\[ \phi_{q_1, \ldots, q_r} : 
\what{\til{\mrm{N}}}{}^{q_1} \mrm{C}(\bsym{U}, \mcal{M}_1)
\times \cdots \times
\what{\til{\mrm{N}}}{}^{q_r} \mrm{C}(\bsym{U}, \mcal{M}_r) \to
\what{\til{\mrm{N}}}{}^q \mrm{C}(\bsym{U}, \mcal{N}) \]
that commute with the embeddings \tup{(\ref{eqn3.1})}.
\end{lem}

\begin{proof}
Direct verification.
\end{proof}

\begin{lem} \label{lem4.9}
Let $\mcal{M}_1, \ldots, \mcal{M}_r, \mcal{N}$ 
be dir-inv $\K_X$-modules, and
$\phi : \prod \mcal{M}_i \to \mcal{N}$
a continuous $\K$-multilinear sheaf homomorphism. 
Then there is an induced homomorphism of complexes of sheaves 
\[ \phi : \what{\til{\mrm{N}}} \mrm{C}(\bsym{U}, \mcal{M}_1)
\otimes \cdots \otimes
\what{\til{\mrm{N}}} \mrm{C}(\bsym{U}, \mcal{M}_r) \to
\what{\til{\mrm{N}}} \mrm{C}(\bsym{U}, \mcal{N}) . \]
\end{lem}

\begin{proof}
Use Lemma \ref{lem3.2}.
\end{proof}

In particular, if $\mcal{M}$ is a dir-inv $\mcal{O}_{X}$-module 
then 
$\what{\til{\mrm{N}}} \mrm{C}(\bsym{U}, \mcal{M})$
is a DG 
$\what{\til{\mrm{N}}} \mrm{C}(\bsym{U}, \mcal{O}_X)$-module.

If $\mcal{M} = \boplus_p \, \mcal{M}^p$ is a graded dir-inv
$\K_X$-module then we define
\[ \what{\til{\mrm{N}}} \mrm{C}(\bsym{U}, \mcal{M})^i :=
\boplus_{p+q=i} \, 
\what{\til{\mrm{N}}}{}^q \mrm{C}(\bsym{U}, \mcal{M}^p)  \]
and
\[ \what{\til{\mrm{N}}} \mrm{C}(\bsym{U}, \mcal{M}) :=
\boplus_{i} \, 
\what{\til{\mrm{N}}} \mrm{C}(\bsym{U}, \mcal{M})^i . \]
Due to Lemma \ref{lem4.9}, if $\mcal{M}$ is a complex in
$\cat{Dir} \cat{Inv} \cat{Mod} \K_X$, then 
$\what{\til{\mrm{N}}} \mrm{C}(\bsym{U}, \mcal{M})$
is also a complex (in $\cat{Mod} \K_X$), and there is a functorial 
homomorphism of complexes
$\mcal{M} \to \what{\til{\mrm{N}}} \mrm{C}(\bsym{U}, \mcal{M})$.

\begin{thm} \label{thm3.1}
Let $X$ be a noetherian topological space, with open covering
$\bsym{U} = \{ U_{(i)} \}_{i = 0}^m$. Let 
$\mcal{M}$ be a bounded below complex in 
$\cat{Dir} \cat{Inv} \cat{Mod} \K_X$, and assume each 
$\mcal{M}^p$ is a complete dir-inv $\K_X$-module. Then:
\begin{enumerate}
\item For any open set $V \subset X$ the homomorphism
\[ \Gamma(V, \int_{\bsym{\Delta}}) :
\Gamma \big( V, \what{\til{\mrm{N}}} \mrm{C}(\bsym{U}, \mcal{M}) 
\big) \to 
\Gamma \big( V, \mrm{N} \mrm{C}(\bsym{U}, \mcal{M}) \big) \]
is a quasi-isomorphism of complexes of $\K$-modules. 
\item There are functorial quasi-isomorphism 
of complexes of $\K_X$-modules
\[ \mcal{M} \to 
\what{\til{\mrm{N}}} \mrm{C}(\bsym{U}, \mcal{M}) 
\xar{\int_{\bsym{\Delta}}}
\mrm{N} \mrm{C}(\bsym{U}, \mcal{M}) . \]
\end{enumerate} 
\end{thm}

\begin{proof}
(1) Lemma \ref{lem3.5} and
Proposition \ref{prop4.4}  imply that for any $p$
the homomorphism of complexes
\[ \Gamma(V, \int_{\bsym{\Delta}}) :
\Gamma \big( V, \what{\til{\mrm{N}}} \mrm{C}(\bsym{U}, \mcal{M}^p) 
\big) \to 
\Gamma \big( V, \mrm{N} \mrm{C}(\bsym{U}, \mcal{M}^p) \big) \]
is a quasi-isomorphism. Now use the standard filtration argument 
(the complexes in question are all bounded below).

\medskip \noindent
(2) From (1) we deduce that
\begin{equation} \label{eqn5.1}
\Gamma(V, \int_{\bsym{\Delta}}) :
\Gamma \big( V, \what{\til{\mrm{N}}} \mrm{C}(\bsym{U}, \mcal{M}) 
\big) \to 
\Gamma \big( V, \mrm{N} \mrm{C}(\bsym{U}, \mcal{M}) \big)
\end{equation}
is a quasi-isomorphism. Hence
\[ \int_{\bsym{\Delta}} :
\what{\til{\mrm{N}}} \mrm{C}(\bsym{U}, \mcal{M}) \to
\mrm{N} \mrm{C}(\bsym{U}, \mcal{M}) \]
is a quasi-isomorphism of complexes of sheaves. 

It is a known fact that 
$\mcal{M}^p \to \mrm{N} \mrm{C}(\bsym{U}, \mcal{M}^p)$
is a quasi-isomorphism of sheaves (see \cite{Ha} Lemma 4.2). 
Again this implies that 
$\mcal{M} \to \mrm{N} \mrm{C}(\bsym{U}, \mcal{M})$
is a quasi-isomorphism. And therefore the homomorphism
$\mcal{M} \to \what{\til{\mrm{N}}} \mrm{C}(\bsym{U}, \mcal{M})$
coming from (\ref{eqn3.3}) is also a quasi-isomorphism.
\end{proof}

Now let us look at a separated noetherian formal scheme 
$\mfrak{X}$. Let $\mcal{I}$ be some defining ideal of $\mfrak{X}$,
and let $X$ be the scheme with structure sheaf 
$\mcal{O}_{X} := \mcal{O}_{\mfrak{X}} / \mcal{I}$. So 
$\mfrak{X}$ and $X$ have the same underlying topological space.
Recall that a dir-coherent $\mcal{O}_{\mfrak{X}}$-module 
is a quasi-coherent $\mcal{O}_{\mfrak{X}}$-module which  
is the union of its coherent submodules.

\begin{cor} \label{cor3.9}
Let $\mfrak{X}$ be a noetherian separated formal scheme over 
$\K$, with defining ideal $\mcal{I}$ and
underlying topological space $X$. Let
$\bsym{U} = \{ U_{(i)} \}_{i = 0}^m$ 
be an affine open covering of $X$. Let 
$\mcal{M}$ be a bounded below complex of sheaves 
of $\K$-modules on $X$. Assume each $\mcal{M}^p$
is a dir-coherent $\mcal{O}_{\mfrak{X}}$-module, and 
the coboundary operators
$\mcal{M}^p \to \mcal{M}^{p+1}$
are continuous for the $\mcal{I}$-adic dir-inv structures
\tup{(}but not necessarily $\mcal{O}_{\mfrak{X}}$-linear\tup{)}.
Then:
\begin{enumerate}
\item The canonical morphism
\[ \Gamma \big( X, \what{\til{\mrm{N}}} \mrm{C}(\bsym{U}, \mcal{M}) 
\big) \to 
\mrm{R} \Gamma \big( X, \what{\til{\mrm{N}}} 
\mrm{C}(\bsym{U}, \mcal{M}) \big) \]
in $\msf{D}(\cat{Mod} \K)$ is an isomorphism. 
\item There is a functorial isomorphism 
\[ \Gamma \big( X, \what{\til{\mrm{N}}} \mrm{C}(\bsym{U}, \mcal{M}) 
\big) \cong \mrm{R} \Gamma (X, \mcal{M}) \]
in $\msf{D}(\cat{Mod} \K)$.
\end{enumerate}
\end{cor}

\begin{proof}
(1) Consider the commutative diagram 
\[ \UseTips \xymatrix @C=10ex @R=5ex {
\Gamma \big( X, \what{\til{\mrm{N}}} \mrm{C}(\bsym{U}, \mcal{M}) 
\big)
\ar[r]^{\Gamma(X, \int_{\bsym{\Delta}})}
\ar[d] 
& \Gamma \big( X, \mrm{N} \mrm{C}(\bsym{U}, \mcal{M}) \big)
\ar[d]
\\
\mrm{R} \Gamma \big( X, \what{\til{\mrm{N}}} 
\mrm{C}(\bsym{U}, \mcal{M}) \big)
\ar[r]^{\mrm{R} \Gamma(X, \int_{\bsym{\Delta}})}
& \mrm{R} \Gamma \big( X, \mrm{N} \mrm{C}(\bsym{U}, \mcal{M}) \big)
} \]
in $\msf{D}(\cat{Mod} \K)$, in which the vertical arrows are the 
canonical morphisms. By part (1) of the theorem (with $V = X$) the 
top arrow is a quasi-isomorphism. And by part (2) the bottom arrow 
is an isomorphism. Hence it is enough to prove that the right 
vertical arrow is an isomorphism. 

Using a filtration argument we may assume that $\mcal{M}$ is a 
single dir-coherent $\mcal{O}_{\mfrak{X}}$-module. 
Now 
$\Gamma \big( X, \mrm{N} \mrm{C}(\bsym{U}, \mcal{M}) \big)$
is the usual \v{C}ech resolution of the sheaf $\mcal{M}$ with 
respect to the covering $\bsym{U}$ (cf.\ equation \ref{eqn3.2})).
So it suffices to prove that for all $q$ and  
$\bsym{i} \in \bsym{\Delta}_q^{m, \mrm{nd}}$ the sheaves 
$g_{\bsym{i} *}\, g_{\bsym{i}}^{-1} \mcal{M}$
are $\Gamma(X, -)$-acyclic.

First let's assume $\mcal{M}$ is a coherent 
$\mcal{O}_{\mfrak{X}}$-module. Let $\mfrak{U}_{\bsym{i}}$
be the open formal subscheme of $\mfrak{X}$ supported on 
$U_{\bsym{i}}$. Then $g_{\bsym{i}}^{-1} \mcal{M}$ is a coherent 
$\mcal{O}_{\mfrak{U}_{\bsym{i}}}$-module, and both 
$g_{\bsym{i}} : \mfrak{U}_{\bsym{i}} \to \mfrak{X}$
and
$\mfrak{U}_{\bsym{i}} \to \opn{Spec} \K$
are affine morphisms. By \cite[Theorem 10.10.2]{EGA-I}  it follows 
that 
$g_{\bsym{i} *}\, g_{\bsym{i}}^{-1} \mcal{M} 
= \mrm{R} g_{\bsym{i} *}\, g_{\bsym{i}}^{-1} \mcal{M}$,
and also
\[ \Gamma(U_{\bsym{i}}, g_{\bsym{i}}^{-1} \mcal{M}) = 
\mrm{R} \Gamma(U_{\bsym{i}}, g_{\bsym{i}}^{-1} \mcal{M}) \cong
\mrm{R} \Gamma(X, \mrm{R} g_{\bsym{i} *} \, g_{\bsym{i}}^{-1}
\mcal{M}) \cong
\mrm{R} \Gamma(X, g_{\bsym{i} *}\, g_{\bsym{i}}^{-1} \mcal{M}) 
. \]
We conclude that 
$\mrm{H}^j (X, g_{\bsym{i} *}\, g_{\bsym{i}}^{-1} \mcal{M}) = 0$
for all $j > 0$. 

In the general case when $\mcal{M}$ is a direct limit of coherent 
$\mcal{O}_{\mfrak{X}}$-modules we still get
$\mrm{H}^j (X, g_{\bsym{i} *}\, g_{\bsym{i}}^{-1} \mcal{M}) = 0$
for all $j > 0$.

\medskip \noindent
(2) By part (2) of the theorem we get a functorial isomorphism 
$\mrm{R} \Gamma (X, \mcal{M}) \cong
\mrm{R} \Gamma \big( X, \what{\til{\mrm{N}}} 
\mrm{C}(\bsym{U}, \mcal{M}) \big)$. 
Now use part (1) above.
\end{proof}

\section{Mixed Resolutions}
\label{sec4}

In this section $\K$ is s field of characteristic $0$
and $X$ is a finite type $\K$-scheme. 

Let us begin be recalling the definition of the sheaf of principal 
parts $\mcal{P}_X$ from \cite{EGA-IV}. Let
$\Delta : X \to X^2 = X \times_{\K} X$ be the diagonal embedding. 
By completing $X^2$ along $\Delta(X)$ we obtain a noetherian 
formal scheme $\mfrak{X}$, and
$\mcal{P}_X := \mcal{O}_{\mfrak{X}}$. 
The two projections
$\mrm{p}_i : X^2 \to X$ give rise to two ring homomorphisms
$\mrm{p}_i^* : \mcal{O}_X \to \mcal{P}_X$. We view 
$\mcal{P}_X$ as a left (resp.\ right) $\mcal{O}_{X}$-module via
$\mrm{p}_1^*$ (resp.\ $\mrm{p}_2^*$).

Recall that a connection $\nabla$ on an $\mcal{O}_{X}$-module 
$\mcal{M}$ is a $\K$-linear sheaf homomorphism
$\nabla : \mcal{M} \to \Omega^{1}_X \otimes_{\mcal{O}_{X}} 
\mcal{M}$
satisfying the Leibniz rule
$\nabla(f m) = \d(f) \otimes m + f \nabla(m)$
for local sections $f \in \mcal{O}_{X}$ and $m \in \mcal{M}$.

\begin{dfn} \label{dfn2.4}
Consider the de Rham differential 
$\mrm{d}_{X^2 / X} : \mcal{O}_{X^2} \to \Omega^1_{X^2 / X}$
relative to the morphism $\mrm{p}_2 : X^2 \to X$.
Since 
$\Omega^1_{X^2 / X} \cong \mrm{p}^*_1\, \Omega^1_{X} = 
\mrm{p}^{-1}_1 \Omega^1_{X} 
\otimes_{\mrm{p}^{-1}_1 \mcal{O}_X} \mcal{O}_{X^2}$
we obtain a $\K$-linear homomorphism
$\mrm{d}_{X^2 / X}  : \mcal{O}_{X^2} \to 
\mrm{p}^*_1 \, \Omega^1_{X}$. 
Passing to the completion along the diagonal $\Delta(X)$ 
we get a connection of $\mcal{O}_X$-modules
\begin{equation} \label{eqn1.1}
\nabla_{\mcal{P}} : \mcal{P}_{X} \to \Omega^1_{X}
\otimes_{\mcal{O}_X} \mcal{P}_{X} 
\end{equation}
called the {\em Grothendieck connection}.
\end{dfn}

Note that the connection 
$\nabla_{\mcal{P}}$ is $\mrm{p}^{-1}_2 \mcal{O}_X$ -linear.
It will be useful to describe $\nabla_{\mcal{P}}$
on the level of rings. Let
$U = \opn{Spec} C \subset X$ be an affine open set. 
Then 
\[ \Gamma(U, \Omega^1_{X} \otimes_{\mcal{O}_X} \mcal{P}_{X}) \cong 
\Omega^{1}_C \otimes_C (\what{C \otimes C}) \cong
\what{\Omega^{1}_C  \otimes C} , \]
the $I$-adic completion, where
$I := \opn{Ker}(C \otimes C \to C)$. 
And 
$\nabla_{\mcal{P}} : \what{C \otimes C} \to
\what{\Omega^{1}_C \otimes C}$
is the completion of 
$\mrm{d} \otimes \bsym{1} : C \otimes C \to 
\Omega^{1}_C \otimes C$.

As usual the connection $\nabla_{\mcal{P}}$ of (\ref{eqn1.1}) 
induces differential operators of left $\mcal{O}_{X}$-modules
\[ \nabla_{\mcal{P}} : \Omega^i_{X} 
\otimes_{\mcal{O}_X} \mcal{P}_{X} 
\to \Omega^{i + 1}_{X} \otimes_{\mcal{O}_X} \mcal{P}_{X} \]
for all $i \geq 0$, by the rule
\begin{equation} \label{eqn1.5}
\nabla_{\mcal{P}}(\alpha \otimes b) = \mrm{d}(\alpha) \otimes b
+ (-1)^i \alpha \wedge \nabla_{\mcal{P}}(b) . 
\end{equation} 

\begin{thm} \label{thm1.5}
Assume $X$ is a smooth $n$-dimensional $\K$-scheme.
Let $\mcal{M}$ be an $\mcal{O}_X$-module. 
Then the sequence of sheaves on $X$
\begin{equation} \label{eqn1.4}
\begin{aligned}
0 \to \mcal{M} & \xar{m \mapsto 1 \otimes m}
\mcal{P}_{X} \otimes_{\mcal{O}_X} \mcal{M}
\xar{\nabla_{\mcal{P}} \otimes \bsym{1}_{\mcal{M}}} 
\Omega^1_{X} \otimes_{\mcal{O}_X} 
\mcal{P}_{X}  \otimes_{\mcal{O}_X} \mcal{M} \\
& \cdots 
\xar{\nabla_{\mcal{P}} \otimes \bsym{1}_{\mcal{M}}} 
\Omega^n_{X} \otimes_{\mcal{O}_X} \mcal{P}_{X} 
\otimes_{\mcal{O}_X} \mcal{M} \to 0
\end{aligned}
\end{equation}
is exact.
\end{thm} 

\begin{proof}
The proof is similar to that of \cite[Theorem 4.5]{Ye1}.
We may restrict to an affine open set 
$U = \opn{Spec} B \subset X$ that 
admits an \'etale coordinate system $\bsym{s} = (s_1, \ldots, s_n)$,
i.e.\ $\K[\bsym{s}] \to B$ is an \'etale ring homomorphism.
It will be convenient to have another copy of $B$, which we call 
$C$; so that
$\Gamma(U, \mcal{P}_X) = \what{B \otimes C}$,
the $I$-adic completion, where
$I := \opn{Ker}(B \otimes C \to B)$. 
We shall  identify $B$ and $C$ with their images inside
$B \otimes C$, and denote the copy of the element 
$s_i$ in $C$ by $r_i$. Letting
$t_i := r_i - s_i \in B \otimes C$ we then have 
$t_i = \til{s}_i = 1 \otimes s_i - s_i \otimes 1$
in our earlier notation. Note that 
$\Omega^{}_{\K[\bsym{s}]} \subset \Omega^{}_B$
is a sub DG algebra, and 
$B \otimes_{\K[\bsym{s}]} \Omega^{}_{\K[\bsym{s}]} \to \Omega^{}_B$
is a bijection. 

By definition
\begin{equation} \label{eqn4.4}
\Gamma(U, \Omega^{}_X \otimes_{\mcal{O}_{X}}
\mcal{P}_X) \cong \Omega^{}_B \otimes_B 
(\what{B \otimes C}) \cong
\what{\Omega^{}_B \otimes C} . 
\end{equation}
The differential $\nabla_{\mcal{P}}$ on the left goes to the 
differential $\mrm{d}_B \otimes \bsym{1}_C$ on the right.
Consider the sub DG algebra
$\Omega^{}_{\K[\bsym{s}]} \otimes C \subset
\Omega^{}_{B} \otimes C$. 
We know that $\K \to \Omega^{}_{\K[\bsym{s}]}$ is a 
quasi-isomorphism; therefore so is
$C \to \Omega^{}_{\K[\bsym{s}]} \otimes C$.

Because $t_i + s_i = r_i \in C$ we see that 
$C[\bsym{s}] = C[\bsym{t}] \subset B \otimes C$.
Therefore we obtain $C$-linear isomorphisms
\[ \Omega^{p}_{\K[\bsym{s}]} \otimes C \cong
\Omega^{p}_{\K[\bsym{s}]} \otimes_{\K[\bsym{s}]} C[\bsym{s}]
= \Omega^{p}_{\K[\bsym{s}]} \otimes_{\K[\bsym{s}]} C[\bsym{t}] . 
\]
So there is a commutative diagram
\begin{equation} \label{eqn2.8}
\UseTips \xymatrix @=4ex {
0 \ar[r] & C \ar[r] \ar[d]
& C[\bsym{t}] \ar[r]^<<<<<{\nabla_{\mcal{P}}} \ar[d]
& \Omega^{1}_{\K[\bsym{s}]} \otimes_{\K[\bsym{s}]} C[\bsym{t}]
\ar[r]^{\nabla_{\mcal{P}} \ \ } \ar@<-0.5ex>[d]
& \ \cdots \ 
\Omega^{n}_{\K[\bsym{s}]} \otimes_{\K[\bsym{s}]} C[\bsym{t}]
\ar[r] \ar@<1ex>[d] & 0 \\
0 \ar[r] & C \ar[r] 
& B \otimes C \ar[r]^>>>>>{\nabla_{\mcal{P}}} & 
\, \Omega^{1}_{B} \otimes_{} C \quad
\ar[r]^{\nabla_{\mcal{P}}} 
& \cdots 
\ \Omega^{n}_{B} \otimes C \quad \ar[r] & 0 }
\end{equation}
of $C$-modules. The top row is exact, and the vertical arrow are 
inclusions. 
Let us introduce a new grading on 
$\Omega^{p}_{\K[\bsym{s}]} \otimes_{\K[\bsym{s}]} C[\bsym{t}]$
as follows: 
$\opn{deg}(s_i) := 1$, $\opn{deg}(t_i) := 1$, 
$\opn{deg}(\mrm{d}(s_i)) := 1$ and $\opn{deg}(c) := 0$
for every nonzero $c \in C$. 
Since $\nabla_{\mcal{P}}(t_i) = -\mrm{d} (s_i)$
we see that $\nabla_{\mcal{P}}$ is homogeneous of degree $0$, and 
thus the top row in (\ref{eqn2.8}) is an exact sequence in the 
category $\cat{GrMod} C$ of graded $C$-modules. Now each term in 
this sequence is a free graded $C$-module, 
and therefore this sequence is split in $\cat{GrMod} C$.

The $\bsym{t}$-adic inv structure on $C[\bsym{t}]$ can be recovered 
from the grading, and this inv structure 
is the same as the $I$-adic inv structure on $B \otimes C$. 
Therefore the completion is
$\Omega^{p}_{\K[\bsym{s}]} \otimes_{\K[\bsym{s}]} C[[\bsym{t}]]
\cong \what{\Omega^p_B \otimes C}$. 
Thus the diagram (\ref{eqn2.8}) is transformed to the commutative 
diagram 
\[ \UseTips \xymatrix @=4ex {
0 \ar[r] & C \ar[r] \ar[d]
& C[[\bsym{t}]] \ar[r]^(0.3){\nabla_{\mcal{P}}} \ar[d]
& \Omega^{1}_{\K[\bsym{s}]} \otimes_{\K[\bsym{s}]} C[[\bsym{t}]]
\ar[r]^(0.45){\nabla_{\mcal{P}}} \ar[d]
& \ \cdots \
\Omega^{n}_{\K[\bsym{s}]} \otimes_{\K[\bsym{s}]} C[[\bsym{t}]]
\ar[r] \ar@<1.2ex>[d] & 0 \\
0 \ar[r] & C \ar[r] 
& \what{B \otimes C} \ar[r]^(0.4){\nabla_{\mcal{P}}} & 
\ \what{\Omega^{1}_{B} \otimes C} \
\ar[r]^(0.45){\nabla_{\mcal{P}}} 
& \ \cdots \ \what{ \Omega^{n}_{B} \otimes C} 
\ \ar[r] & 0 
} \]
in which the top tow is continuously $C$-linearly split, and the 
vertical arrows are bijections. 
Hence the bottom row is split exact. Comparing this to 
(\ref{eqn4.4}) we conclude that the sequence of right 
$\mcal{O}_U$-modules
\[ 0 \to \mcal{O}_U \xar{\mrm{p}^*_2} \mcal{P}_{X}|_U 
\xar{\nabla_{\mcal{P}}} \bigl( \Omega^1_{X} \otimes_{\mcal{O}_X}
\mcal{P}_{X} \bigr)|_U \xar{\nabla_{\mcal{P}}} \cdots 
\bigl( \Omega^n_{X} \otimes_{\mcal{O}_X} \mcal{P}_{X}
\bigr)|_U \to 0 \]
is split exact. 

Therefore it follows that for any $\mcal{O}_X$-module 
$\mcal{M}$ the sequence (\ref{eqn1.4}), when restricted to $U$, is 
split exact.
\end{proof}

Let us now fix an affine open covering
$\bsym{U} = \{ U_{(0)}, \ldots, U_{(m)} \}$ of $X$. 

Let 
$\mcal{I}_X = \opn{Ker}(\mcal{P}_X \to \mcal{O}_X)$.
This is a defining ideal of the noetherian formal scheme
$(\mfrak{X}, \mcal{O}_{\mfrak{X}}) := (X, \mcal{P}_X)$.
So $\mcal{P}_X$ is an inv module over itself
with the $\mcal{I}_X$-adic inv structure. 
Given quasi-coherent $\mcal{O}_{X}$-modules
$\mcal{M}$ and $\mcal{N}$, the tensor product
$\mcal{N} \otimes_{\mcal{O}_{X}} \mcal{P}_X \otimes_{\mcal{O}_{X}}
\mcal{M}$
is a dir-coherent $\mcal{P}_X$-module, and so it has the 
$\mcal{I}_{X}$-adic dir-inv structure. See Example \ref{exa2.1}.
In particular 
\[ \Omega_{X} \otimes_{\mcal{O}_{X}} \mcal{P}_X 
\otimes_{\mcal{O}_{X}} \mcal{M} = 
\boplus_{p \geq 0}
\Omega^{p}_{X} \otimes_{\mcal{O}_{X}} \mcal{P}_X 
\otimes_{\mcal{O}_{X}} \mcal{M} \]
becomes a dir-inv $\K_X$-module.

\begin{lem}
$\Omega_{X} \otimes_{\mcal{O}_{X}} \mcal{P}_X 
\otimes_{\mcal{O}_{X}} \mcal{M}$
is a DG $\Omega_{X}$-module in 
$\cat{Dir} \cat{Inv} \cat{Mod} \K_X$, with differential
$\nabla_{\mcal{P}} \otimes \bsym{1}_{\mcal{M}}$.
\end{lem}

\begin{proof}
Because 
$\nabla_{\mcal{P}} \otimes \bsym{1}_{\mcal{M}}$
is a differential operator of $\mcal{P}_X$-modules, it is 
continuous for the $\mcal{I}_X$-adic dir-inv structure. See 
\cite[Proposition 2.3]{Ye2}. 
\end{proof}

Henceforth we will write $\nabla_{\mcal{P}}$ instead of
$\nabla_{\mcal{P}} \otimes \bsym{1}_{\mcal{M}}$.

\begin{dfn} \label{dfn4.1}
Let $\mcal{M}$ be a quasi-coherent $\mcal{O}_X$-module. 
For any $p, q \in \mbb{N}$ define
\[ \opn{Mix}_{\bsym{U}}^{p, q}(\mcal{M}) := 
\what{\til{\mrm{N}}}{}^q \mrm{C}(\bsym{U}, \Omega^{p}_{X}
\otimes_{\mcal{O}_X} \mcal{P}_X \otimes_{\mcal{O}_X}
\mcal{M}) . \]
The Grothendieck connection
\[ \nabla_{\mcal{P}} : 
\Omega^{p}_{X} \otimes_{\mcal{O}_X} \mcal{P}_X 
\otimes_{\mcal{O}_X} \mcal{M} \to
\Omega^{p+1}_{X} \otimes_{\mcal{O}_X} \mcal{P}_X
\otimes_{\mcal{O}_X} \mcal{M} \]
induces a homomorphism of sheaves
\[ \nabla_{\mcal{P}} : \opn{Mix}_{\bsym{U}}^{p, q}(\mcal{M}) \to 
\opn{Mix}_{\bsym{U}}^{p+1, q}(\mcal{M}) . \]
We also have 
$\partial : \opn{Mix}_{\bsym{U}}^{p, q}(\mcal{M}) \to 
\opn{Mix}_{\bsym{U}}^{p, q+1}(\mcal{M})$.
Define 
\[ \opn{Mix}_{\bsym{U}}^{i}(\mcal{M}) := 
\boplus_{p + q = i}\, \opn{Mix}_{\bsym{U}}^{p, q}(\mcal{M}) , \]
\[ \opn{Mix}^{}_{\bsym{U}}(\mcal{M}) := 
\boplus_{i}\, \opn{Mix}_{\bsym{U}}^{i}(\mcal{M}) \]
and
\begin{equation} \label{eqn4.3}
\d_{\mrm{mix}} := \partial + (-1)^q \nabla_{\mcal{P}} : 
\opn{Mix}_{\bsym{U}}^{p, q}(\mcal{M}) \to 
\opn{Mix}_{\bsym{U}}^{p+1, q} \oplus
\opn{Mix}_{\bsym{U}}^{p, q+1}(\mcal{M}) . 
\end{equation}
The complex 
$\big( \opn{Mix}_{\bsym{U}}(\mcal{M}), \d_{\mrm{mix}} \big)$
is called the {\em mixed resolution of $\mcal{M}$}.
\end{dfn}

There are functorial embeddings of sheaves
\begin{equation} \label{eqn5.2}
\mcal{M} \subset \mcal{P}_X \otimes_{\mcal{O}_{X}} \mcal{M} 
\subset
\what{\til{\mrm{N}}}{}^0 \mrm{C}(\bsym{U}, \Omega^{0}_{X}
\otimes_{\mcal{O}_X} \mcal{P}_X \otimes_{\mcal{O}_X}
\mcal{M}) =
\opn{Mix}_{\bsym{U}}^{0, 0}(\mcal{M}) 
\end{equation}
and
\begin{equation} \label{eqn5.3}
\opn{Mix}^{p,q}_{\bsym{U}}(\mcal{M})
\subset 
{\displaystyle \prod_{l \in \mbb{N}}}\ \
{\displaystyle \prod_{\bsym{i} \in 
\bsym{\Delta}_l^{m}}} \,
g_{\bsym{i} *}\, g_{\bsym{i}}^{-1}\,
\big( \Omega^{q}(\bsym{\Delta}^l_{\K})
\, \what{\otimes} \, (\Omega^{p}_X \otimes_{\mcal{O}_{X}}
\mcal{P}_X \otimes_{\mcal{O}_{X}} \mcal{M}) \big) ;
\end{equation}
see Lemma \ref{lem3.5}.

\begin{prop} \label{prop5.2}  
\begin{enumerate}
\item 
$\opn{Mix}^{}_{\bsym{U}}(\mcal{O}_X)$ 
is a sheaf of super-commutative associative unital
DG $\K$-algebras. There are two 
$\K$-algebra homomorphisms
$\mrm{p}_1^*, \mrm{p}_2^* : \mcal{O}_X \to 
\opn{Mix}^0_{\bsym{U}}(\mcal{O}_X)$.
\item Let $\mcal{M}$ be a quasi-coherent $\mcal{O}_X$-module. Then 
$\opn{Mix}^{}_{\bsym{U}}(\mcal{M})$
is a left DG $\opn{Mix}^{}_{\bsym{U}}(\mcal{O}_X)$-module. 
\item If $\mcal{M}$ is a locally free $\mcal{O}_X$-module of 
finite rank then the multiplication map
\[ \opn{Mix}^{}_{\bsym{U}}(\mcal{O}_X) 
\otimes_{\mcal{O}_X} \mcal{M}
\to \opn{Mix}^{}_{\bsym{U}}(\mcal{M}) \]
is an isomorphism.
\end{enumerate}
\end{prop}

\begin{proof} 
By by Lemmas \ref{lem3.5} and \ref{lem4.9}.
\end{proof}

Note that 
$\d_{\mrm{mix}} \circ \mrm{p}_2^* : \mcal{O}_{X} \to
\opn{Mix}^{}_{\bsym{U}}(\mcal{O}_X)$
is zero, but $\d_{\mrm{mix}} \circ \mrm{p}_1^* \neq 0$.

\begin{prop} \label{prop5.1}
Let $\mcal{M}_1, \ldots, \mcal{M}_r, \mcal{N}$ be quasi-coherent 
$\mcal{O}_{X}$-modules. Suppose  
\[ \phi : \prod_{i=1}^r 
(\Omega_X \otimes_{\mcal{O}_{X}} \mcal{P}_X \otimes_{\mcal{O}_{X}}
\mcal{M}_i) \to 
\Omega_X \otimes_{\mcal{O}_{X}} \mcal{P}_X \otimes_{\mcal{O}_{X}}
\mcal{N} \]
is a continuous $\Omega_X$-multilinear sheaf morphism of degree $d$. 
Then there is a unique $\K$-multilinear sheaf morphism of 
degree $d$ 
\[ \what{\til{\mrm{N}}} \mrm{C}(\bsym{U}, \phi) :
\opn{Mix}_{\bsym{U}}(\mcal{M}_1) \times \cdots
\times \opn{Mix}_{\bsym{U}}(\mcal{M}_r) \to
\opn{Mix}_{\bsym{U}}(\mcal{N}) \]
which is compatible with 
$\phi$ via the embedding \tup{(\ref{eqn5.3})}.
\end{prop}

\begin{proof}
This is an immediate consequence of Lemma \ref{lem4.9}.
\end{proof}

Suppose we are given
$\mcal{M} \in \msf{C}^{+}(\cat{QCoh} \mcal{O}_X)$. 
Define
\[ \opn{Mix}_{\bsym{U}}(\mcal{M})^i :=
\bigoplus_{p+q=i} \opn{Mix}^q_{\bsym{U}}(\mcal{M}^p) \]
with differential
\[ \d_{\mrm{mix}} + (-1)^q \d_{\mcal{M}} :
\opn{Mix}^q_{\bsym{U}}(\mcal{M}^p) \to
\opn{Mix}^{q+1}_{\bsym{U}}(\mcal{M}^p) \oplus
\opn{Mix}^q_{\bsym{U}}(\mcal{M}^{p+1}) . \]

\begin{thm} \label{thm4.1} 
Let $X$ be a smooth separated $\K$-scheme, and let 
$\bsym{U} = \linebreak \{ U_{(0)}, \ldots, U_{(m)} \}$ be an affine
open covering of $X$.
\begin{enumerate}
\item There is a functorial quasi-isomorphism
$\mcal{M} \to \opn{Mix}^{}_{\bsym{U}}(\mcal{M})$
for $\mcal{M} \in \linebreak \msf{C}^{+}(\cat{QCoh} \mcal{O}_X)$.
\item Given 
$\mcal{M} \in \msf{C}^{+}(\cat{QCoh} \mcal{O}_X)$,
the canonical morphism 
$\Gamma \big( X, \opn{Mix}^{}_{\bsym{U}}(\mcal{M}) \big) \to
\mrm{R} \Gamma \bigl( X, \opn{Mix}^{}_{\bsym{U}}(\mcal{M}) \big)$
in $\msf{D}(\cat{Mod} \K)$ is an isomorphism.
\item The quasi-isomorphism in part \tup{(1)} induces a 
functorial isomorphism \linebreak
$\Gamma \big( X, \opn{Mix}^{}_{\bsym{U}}(\mcal{M}) \big) \cong
\mrm{R} \Gamma(X, \mcal{M})$ 
in $\msf{D}(\cat{Mod} \K)$.
\end{enumerate}
\end{thm}

\begin{proof}
(1) Write
$\mcal{N} := \Omega_X \otimes_{\mcal{O}_{X}} 
\mcal{P}_X \otimes_{\mcal{O}_{X}} \mcal{M}$. 
A filtration argument and Theorem \ref{thm1.5} show that the 
inclusion $\mcal{M} \to \mcal{N}$ is a quasi-isomorphism. 
Next we view $\mcal{N}$ as a bounded below 
complex in $\cat{Dir} \cat{Inv} \cat{Mod} \K_X$. 
By Theorem \ref{thm3.1}(2)
we have a quasi-isomorphism
$\mcal{N} \to 
\what{\til{\mrm{N}}}{}^q \mrm{C}(\bsym{U}, \mcal{N}) = 
\opn{Mix}^{}_{\bsym{U}}(\mcal{M})$. 

\medskip \noindent
(2) This is due to Corollary \ref{cor3.9}(1), applied to the 
formal scheme $(X, \mcal{P}_X)$ and the complex $\mcal{N}$ of 
dir-coherent $\mcal{P}_X$-modules defined above.

\medskip \noindent
(3) This assertion is an immediate consequence of parts (1) and (2). 
\end{proof}

\begin{cor} \label{cor4.1}
In the situation of the theorem, suppose 
$\mcal{M}, \mcal{N} \in \msf{C}^{+}(\cat{QCoh} \mcal{O}_X)$
and 
$\phi : \opn{Mix}^{}_{\bsym{U}}(\mcal{M}) \to 
\opn{Mix}^{}_{\bsym{U}}(\mcal{N})$
is a $\K$-linear quasi-isomorphism. Then
\[ \Gamma(X, \phi) : 
\Gamma \big( X, \opn{Mix}^{}_{\bsym{U}}(\mcal{M}) \big) \to
\Gamma \big( X, \opn{Mix}^{}_{\bsym{U}}(\mcal{N}) \big) \]
is a quasi-isomorphism.
\end{cor}

\begin{proof}
Consider the commutative diagram
\[ \UseTips \xymatrix @C=8ex @R=5ex {
\Gamma \big( X, \opn{Mix}^{}_{\bsym{U}}(\mcal{M}) \big)
\ar[r]^{\Gamma(X, \phi)}
\ar[d]
& \Gamma \big( X, \opn{Mix}^{}_{\bsym{U}}(\mcal{N}) \big)
\ar[d] \\
\mrm{R} \Gamma \big( X, \opn{Mix}^{}_{\bsym{U}}(\mcal{M}) \big)
\ar[r]^{\mrm{R} \Gamma(X, \phi)}
& \mrm{R} \Gamma \big( X, \opn{Mix}^{}_{\bsym{U}}(\mcal{N}) \big)
} \]
in $\msf{D}(\cat{Mod} \K)$. By part (2) of the theorem the 
vertical arrows are isomorphisms. Since $\phi$ is an isomorphism 
in $\msf{D}(\cat{Mod} \K_X)$ it follows that the bottom arrow is 
an isomorphism.
\end{proof}

Given a quasi-coherent $\mcal{O}_X$-module $\mcal{M}$
and an integer $i$ define 
\[ \mrm{G}^i \opn{Mix}^{}_{\bsym{U}}(\mcal{M}) :=
\bigoplus_{q \geq i}\, \opn{Mix}^{q}_{\bsym{U}}(\mcal{M}) . \]
Then
$\{ \mrm{G}^i \opn{Mix}^{}_{\bsym{U}}(\mcal{M}) \}_{i \in \mbb{Z}}$
is a descending filtration of
$\opn{Mix}^{}_{\bsym{U}}(\mcal{M})$
by subcomplexes, satisfying
$\mrm{G}^i \opn{Mix}^{}_{\bsym{U}}(\mcal{M}) = 
\opn{Mix}^{}_{\bsym{U}}(\mcal{M})$
for $i \ll 0$ and
$\bigcap_i \mrm{G}^i \opn{Mix}^{}_{\bsym{U}}(\mcal{M}) = 0$.
For any $i$ define 
\[ \opn{gr}^i_{\mrm{G}} \opn{Mix}^{}_{\bsym{U}}(\mcal{M}) :=
\mrm{G}^i \opn{Mix}^{}_{\bsym{U}}(\mcal{M})\, /\, 
\mrm{G}^{i+1} \opn{Mix}^{}_{\bsym{U}}(\mcal{M}) . \]
The functor
\[ \opn{gr}^i_{\mrm{G}} \opn{Mix}^{}_{\bsym{U}}
 : \cat{QCoh} \mcal{O}_X \to
\cat{Mod} \K_X \]
is additive, but we do not know whether it is exact. The next 
theorem asserts this in a very special case.

Consider the sheaves of DG Lie algebras 
$\mcal{T}^{}_{\mrm{poly}, X}$ and $\mcal{D}^{}_{\mrm{poly}, X}$
as complexes of quasi-coherent $\mcal{O}_X$-modules
(cf.\ \cite[Proposition 3.18]{Ye3}).
According to \cite[Theorem 0.4]{Ye1} there is a quasi-isomorphism
\[ \mcal{U}_1 : \mcal{T}^{}_{\mrm{poly}, X} \to
\mcal{D}^{}_{\mrm{poly}, X} . \]

\begin{thm} \label{thm4.2}
For any $i$ the homomorphism of complexes
\[ \opn{gr}^i_{\mrm{G}} \opn{Mix}^{}_{\bsym{U}}(\mcal{U}_1) :
\opn{gr}^i_{\mrm{G}} \opn{Mix}^{}_{\bsym{U}} 
(\mcal{T}^{}_{\mrm{poly}, X}) \to
\opn{gr}^i_{\mrm{G}} \opn{Mix}^{}_{\bsym{U}} 
(\mcal{D}^{}_{\mrm{poly}, X}) \]
is a quasi-isomorphism. 
\end{thm}

\begin{proof}
Given a point $x \in X$ choose an affine open neighborhood $V$ of 
$x$ which admits an \'etale morphism $V \to \mbf{A}^n_{\K}$. 
By \cite[Theorem 4.11]{Ye2} the map of complexes
\[ \mcal{U}_1|_V : \mcal{T}^{}_{\mrm{poly}, X}|_V \to
\mcal{D}^{}_{\mrm{poly}, X}|_V \]
is a homotopy equivalence in
$\msf{C}^{+}(\cat{QCoh} \mcal{O}_V)$.
Since $\opn{gr}^i_{\mrm{G}} \opn{Mix}^{}_{\bsym{U}}$ is an 
additive functor we see that 
$\opn{gr}^i_{\mrm{G}} \opn{Mix}^{}_{\bsym{U}}(\mcal{U}_1)|_V$
is a quasi-isomorphism.
\end{proof}

\begin{rem}
We know very little about the structure of the sheaves
$\what{\til{\mrm{N}}}{}^{q} \mrm{C}(\bsym{U}, \mcal{M})$,
even when $\mcal{M} = \mcal{O}_X$. Cf.\ \cite{HS}.
\end{rem}

\section{Simplicial Sections}
\label{sec5}

Let $X$ be a $\K$-scheme, and let
$X = \bigcup_{i = 0}^m U_{(i)}$ be an open covering, with 
inclusions $g_{(i)} : U_{(i)} \to X$. We denote this covering by 
$\bsym{U}$. For any multi-index
$\bsym{i} = (i_0, \ldots, i_q) \in \bsym{\Delta}^m_q$
we write $U_{\bsym{i}} := \bigcap_{j = 0}^q U_{(i_j)}$,
and we define the scheme
$U_{q} := \coprod_{\bsym{i} \in \bsym{\Delta}^m_q}
U_{\bsym{i}}.$
Given $\alpha \in \bsym{\Delta}_p^q$ and 
$\bsym{i} \in \bsym{\Delta}^m_q$ there is an inclusion of open 
sets
$\alpha_* : U_{\bsym{i}} \to U_{\alpha_*(\bsym{i})}$. 
These patch to a morphism of schemes
$\alpha_* : U_q \to U_p$,
making $\{ U_q \}_{q \in \mbb{N}}$ into a simplicial 
scheme. The inclusions $g_{(i)} : U_{(i)} \to X$ induce inclusions
$g_{\bsym{i}} : U_{\bsym{i}} \to X$
and morphisms 
$g_q : U_q \to X$; and one has the relations
$g_p \circ \alpha_* = g_q$
for any $\alpha \in \bsym{\Delta}_p^q$.

\begin{dfn} \label{dfn3.2}
Let $\pi : Z \to X$ be a morphism of $\K$-schemes.
A {\em simplicial section of 
$\pi$} based on the covering $\bsym{U}$ is a sequence of morphisms
\[ \bsym{\sigma} = \{ \sigma_q : 
\bsym{\Delta}^q_{\K} \times U_q \to Z \}_{q \in \mbb{N}}  \]
satisfying the following conditions.
\begin{enumerate}
\rmitem{i} For any $q$ the diagram
\[ \begin{CD}
\bsym{\Delta}^q_{\K} \times U_q @>{\sigma_q}>>
Z \\
@V{\mrm{p}_2}VV @V{\pi}VV \\
U_q @>{g_q}>> X 
\end{CD} \]
is commutative.
\rmitem{ii} For any $\alpha \in \bsym{\Delta}_p^q$
the diagram
\[ \UseTips \xymatrix @=4.5ex {
& \bsym{\Delta}^{p}_{\K} \times U_p
\ar[dr]^{\sigma_{p}} 
\\
\bsym{\Delta}^{p}_{\K} \times U_q
\ar[ur]^{\bsym{1} \times \alpha^*} 
\ar[dr]_{\alpha_* \times \bsym{1}} 
& & Z \\
& \bsym{\Delta}^{q}_{\K} \times U_q
\ar[ur]_{\sigma_{q}} 
} \]
is commutative.
\end{enumerate}
\end{dfn}

\begin{figure}
\includegraphics{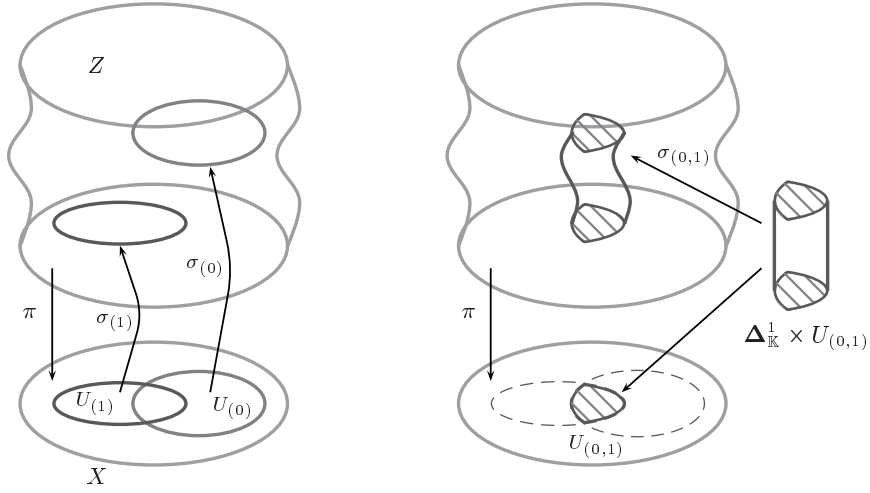}
\caption{An illustration of a simplicial section
$\bsym{\sigma}$ based on an open 
covering $\bsym{U} = \{ U_{(i)} \}$.
On the left we see two components of $\bsym{\sigma}$ in dimension 
$q = 0$; and on the right we see one component in dimension 
$q = 1$.}
\label{fig1}
\end{figure}

Given a multi-index $\bsym{i} \in \bsym{\Delta}^m_q$
we denote by $\sigma_{\bsym{i}}$ the restriction of 
$\sigma_{q}$ to 
$\bsym{\Delta}^q_{\K} \times U_{\bsym{i}}$.
See Figure \ref{fig1} for an illustration. 

As explained in the introduction, simplicial sections arise 
naturally in several contexts, including deformation quantization. 

Let $A$ be an associative unital super-commutative DG $\K$-algebra. 
Consider homogeneous $A$-multilinear functions
$\phi : M_1 \times \cdots \times M_r \to N$,
where $M_1, \ldots, M_r, N$ are DG $A$-modules. There is an 
operation of composition for such functions: given  
functions $\psi_i : \prod_j L_{i,j} \to M_i$ 
the composition is
$\phi \circ (\psi_1 \times \cdots \times \psi_r)
: \prod_{i,j} L_{i,j} \to N$.
There is also a summation operation: if 
$\phi_j : \prod_i M_i \to N$ are homogeneous of equal degree then 
so is their sum $\sum_j \phi_j$.
Finally, let 
$\d : \prod_i M_i \to \prod_i M_i$
be the function 
\[ \d(m_1, \ldots, m_r) :=
\sum_{i=1}^r \pm \bigl( m_1, \ldots, \d(m_i), \ldots, m_r \bigr) \]
with Koszul signs. All the above can of course be sheafified, 
i.e.\ $\mcal{A}$ is a sheaf of DG algebras on a scheme $Z$ etc.

As before let $\pi : Z \to X$ be a morphism if $\K$-schemes,
and let $\bsym{U} = \{ U_{(i)} \}$ be an open covering of $X$. 
Suppose $\bsym{\sigma}$ is a simplicial section of $\pi$ based on 
$\bsym{U}$. We 
consider $\Omega^p_X$ as a discrete inv $\K_X$-module, and 
$\Omega_X = \boplus_{p \geq 0} \Omega^p_X$ has the $\boplus$ 
dir-inv structure. Likewise for 
$\Omega_Z = \boplus_{p \geq 0} \Omega^p_Z$.

Suppose $\mcal{M}$ is a quasi-coherent $\mcal{O}_{X}$-module. 
Then, as explained in Section \ref{sec4}, 
$\Omega_Z \hatotimes{\mcal{O}_{Z}} \pi^{\what{*}}\, 
(\mcal{P}_X \otimes_{\mcal{O}_X} \mcal{M})$
is a DG $\Omega_Z$-module on $Z$, with the Grothendieck connection 
$\nabla_{\mcal{P}}$. And
$\opn{Mix}_{\bsym{U}}^{}(\mcal{M})$ is a DG 
$\opn{Mix}_{\bsym{U}}^{}(\mcal{O}_X)$-module on $X$,
with differential $\d_{\mrm{mix}}$.

\begin{thm} \label{thm3.3}
Let $\pi : Z \to X$ be a morphism of schemes, and suppose
$\bsym{\sigma}$ is a simplicial section of $\pi$ based on an open 
covering $\bsym{U}$ of $X$. 
Let $\mcal{M}_1, \ldots, \mcal{M}_{r}, \mcal{N}$ 
be quasi-coherent $\mcal{O}_X$-modules, and let
\[ \phi  : \prod_{i=1}^r
\big( \Omega^{}_{Z} \, 
\what{\otimes}_{\mcal{O}_{Z}} \,
\pi^{\what{*}}\, 
(\mcal{P}_X \otimes_{\mcal{O}_X} \mcal{M}_i) \big) 
\to \Omega^{}_{Z} \, 
\what{\otimes}_{\mcal{O}_{Z}} \,
\pi^{\what{*}}\, 
(\mcal{P}_X \otimes_{\mcal{O}_X} \mcal{N}) \] 
be a continuous $\Omega_{Z}$-multilinear sheaf
morphism on $Z$ of degree $k$. 
Then there is an induced 
$\opn{Mix}_{\bsym{U}}^{}(\mcal{O}_X)$-multilinear 
sheaf morphism of degree $k$
\[ \bsym{\sigma}^*(\phi) :
\opn{Mix}_{\bsym{U}}^{}(\mcal{M}_1)
\times \cdots \times
\opn{Mix}_{\bsym{U}}^{}(\mcal{M}_r)
\to
\opn{Mix}_{\bsym{U}}^{}(\mcal{N}) \] 
on $X$ with the following properties:
\begin{enumerate}
\rmitem{i} The assignment $\phi \mapsto \bsym{\sigma}^*(\phi)$ 
respects the operations of composition and summation.
\rmitem{ii} If $\phi = \pi^{\what{*}}(\phi_0)$ for some continuous 
$\Omega_X$-multilinear morphism
\[ \phi_0  : \prod_{i=1}^r
\big( \Omega^{}_{X} \otimes_{\mcal{O}_{X}} 
\mcal{P}_X \otimes_{\mcal{O}_X} \mcal{M}_i \big) 
\to \Omega^{}_{X} \otimes_{\mcal{O}_{X}}
\mcal{P}_X \otimes_{\mcal{O}_X} \mcal{N} \] 
then 
$\bsym{\sigma}^*(\phi) = 
\what{\til{\mrm{N}}} \mrm{C}(\bsym{U}, \phi_0)$.
\rmitem{iii} Assume that 
\[ \nabla_{\mcal{P}} \circ \phi - (-1)^k \phi \circ \nabla_{\mcal{P}}
= \psi \]
for some continuous $\Omega_{Z}$-multilinear sheaf
morphism 
\[ \psi  : \prod_{i=1}^r
\big( \Omega^{}_{Z} \, 
\what{\otimes}_{\mcal{O}_{Z}} \,
\pi^{\what{*}}\, 
(\mcal{P}_X \otimes_{\mcal{O}_X} \mcal{M}_i) \big) 
\to \Omega^{}_{Z} \, 
\what{\otimes}_{\mcal{O}_{Z}} \,
\pi^{\what{*}}\, 
(\mcal{P}_X \otimes_{\mcal{O}_X} \mcal{N}) \] 
of degree $k+1$. Then
\[ \d_{\mrm{mix}} \circ \bsym{\sigma}^*(\phi) - (-1)^k 
\bsym{\sigma}^*(\phi) \circ \d_{\mrm{mix}} = 
\bsym{\sigma}^*(\psi) . \]
\end{enumerate}
\end{thm}

Before the proof we need an auxiliary result. 

\begin{lem} \label{lem6.6}
Let $A$ and $B$ be complete DG algebras in 
$\cat{Dir} \cat{Inv} \cat{Mod} \K$,
and let $f^* : A \to B$ be a continuous DG algebra homomorphism.
To any DG $A$-module $M$ in $\cat{Dir} \cat{Inv} \cat{Mod} \K$
we assign the DG $B$-module $f^* M := B \hatotimes{A} M$. 
Then to any continuous $A$-multilinear function
$\phi : \prod_i M_i \to N$ we can assign a continuous 
$B$-multilinear function
$f^*(\phi) : \prod_i f^*(M_i) \to f^*(N)$.
This assignment is functorial in $f^*$, and 
respects the operations of composition and summation.
If $\phi$ and $\psi$ are such continuous $A$-multilinear 
functions, homogeneous of degrees $k$ and $k+1$ respectively and
satisfying 
\[ \d \circ \phi - (-1)^k \phi \circ \d = \psi , \]
then
\[ \d \circ f^*(\phi) - (-1)^k f^*(\phi) \circ \d = 
f^*(\psi) . \]
\end{lem}

\begin{proof}
This is all straightforward, except perhaps the last assertion. 
For that we make the calculations. By continuity and 
multilinearity it suffices to show that
\[ \bigl( \d \circ f^*(\phi) \bigr)(\beta)
 - (-1)^k \bigl( f^*(\phi) \circ \d \bigr)(\beta) = 
f^*(\psi) (\beta) \]
for $\beta = (\beta_1, \ldots, \beta_r)$,
with 
$\beta_i = b_i \otimes m_i$, $b_i \in B^{p_i}$
and $m_i \in M^{q_i}$. Then
\[ \begin{aligned}
\bigl( \d \circ f^*(\phi) \bigr)(\beta) & = 
\d \bigl( \pm b_1 \ldots b_r \cdot \phi(m_1, \ldots, m_r) \bigr) \\
& = \pm \d(b_1 \cdots b_r) \cdot \phi(m_1, \ldots, m_r)
\pm b_1 \cdots b_r \cdot \d \bigl( \phi(m_1, \ldots, m_r) \bigr)
\end{aligned} \]
with Koszul signs. Since
\[ \d(\beta_i) = \d(b_i) \otimes m_i \pm b_i \otimes \d(m_i) \]
we also have
\[ \begin{aligned}
\bigl( f^*(\phi) \circ \d \bigr)(\beta) & =
\sum_i \pm f^*(\phi) \bigl( \beta_1, \ldots, \d(\beta_i),
\ldots, \beta_r) \\
& = \sum_i \Bigl( \pm b_1 \cdots \d(b_i) \cdots b_r \cdot
\phi(m_1, \ldots, m_r) \\
& \qquad \pm b_1 \ldots b_r \cdot \phi \bigl( m_1, \ldots, 
\d(m_i) \cdots m_r \bigr) \Bigr) \\
& = \pm \d(b_1 \cdots b_r) \cdot \phi(m_1, \ldots, m_r)
\pm b_1 \cdots b_r \cdot \phi \bigl( \d(m_1, \ldots, m_r) \bigr) .
\end{aligned} \]
Finally 
\[ f^*(\psi)(\beta) = 
\pm b_1 \cdots b_r \cdot \psi (m_1, \ldots, m_r) , \]
and the signs all match up.
\end{proof}

\begin{proof}[Proof of the theorem]
For a sequence of indices 
$\bsym{i} = (i_0, \ldots, i_l) \in \bsym{\Delta}_l^m$ let 
us introduce the abbreviation 
$Y_{\bsym{i}} := \bsym{\Delta}^l_{\K} \times U_{\bsym{i}}$, and 
let 
$\mrm{p}_2 : Y_{\bsym{i}} \to U_{\bsym{i}}$
be the projection. The simplicial section $\bsym{\sigma}$ 
restricts to a morphism
$\sigma_{\bsym{i}} : Y_{\bsym{i}} \to Z$.

By Lemma \ref{lem6.6}, applied with respect to the DG algebra 
homomorphism 
$\sigma_{\bsym{i}}^* : \sigma_{\bsym{i}}^{-1} \Omega_{Z}
\to  \Omega^{}_{Y_{\bsym{i}}}$,
there is an induced continuous
$\Omega_{Y_{\bsym{i}}}$-multilinear morphism 
\[ \begin{aligned}
& \sigma_{\bsym{i}}^*(\phi) : \prod_{j=1}^r
\Big( \Omega^{}_{Y_{\bsym{i}}} \, 
\what{\otimes}_{\sigma_{\bsym{i}}^{-1} \Omega_{Z}} 
\sigma_{\bsym{i}}^{-1}
\big( \Omega_Z \what{\otimes}_{\mcal{O}_{Z}} \,
\pi^{\what{*}}\, 
(\mcal{P}_X \otimes_{\mcal{O}_X} \mcal{M}_j) \big) \Big) \\
& \qquad \to \ 
\Omega^{}_{Y_{\bsym{i}}} \, 
\what{\otimes}_{\sigma_{\bsym{i}}^{-1} \Omega_{Z}} 
\sigma_{\bsym{i}}^{-1}
\big( \Omega_Z \what{\otimes}_{\mcal{O}_{Z}} \,
\pi^{\what{*}}\, 
(\mcal{P}_X \otimes_{\mcal{O}_X} \mcal{N}) \big) 
\end{aligned} \] 

Now for any quasi-coherent $\mcal{O}_X$-module $\mcal{M}$ 
we have an isomorphism of dir-inv 
DG $\Omega^{}_{Y_{\bsym{i}}}$-modules 
\[ \Omega^{}_{Y_{\bsym{i}}} \, 
\what{\otimes}_{\sigma_{\bsym{i}}^{-1} \Omega_{Z}} 
\sigma_{\bsym{i}}^{-1}
\big( \Omega_Z \what{\otimes}_{\mcal{O}_{Z}} \,
\pi^{\what{*}}\, 
(\mcal{P}_X \otimes_{\mcal{O}_X} \mcal{M}) \big)
\cong
\Omega^{}_{Y_{\bsym{i}}} \what{\otimes}_{\mcal{O}_{Y_{\bsym{i}}}}
\mrm{p}_{2}^* (\mcal{P}_X \otimes_{\mcal{O}_X} \mcal{M}) . \]
Under the DG algebra isomorphism
$\mrm{p}_{2 *} \Omega^{}_{Y_{\bsym{i}}} \cong
\Omega(\bsym{\Delta}^l_{\K}) \otimes \Omega_{U_{\bsym{i}}}$
there is a dir-inv DG module isomorphism
\[ \mrm{p}_{2 *} \bigl( 
\Omega^{}_{Y_{\bsym{i}}} \what{\otimes}_{\mcal{O}_{Y_{\bsym{i}}}}
\mrm{p}_{2}^* (\mcal{P}_X \otimes_{\mcal{O}_X} \mcal{M}) \bigr) 
\cong \Omega(\bsym{\Delta}^l_{\K}) \hatotimes{}
(\Omega_X \otimes_{\mcal{O}_X}
\mcal{P}_X \otimes_{\mcal{O}_X} \mcal{M})|_{U_{\bsym{i}}} . \]
Thus we obtain a family of morphisms
\[ \begin{aligned}
& \sigma_{\bsym{i}}^*(\phi) : \prod_{j=1}^r
\Big( \Omega(\bsym{\Delta}^l_{\K}) \hatotimes{}
(\Omega_X \otimes_{\mcal{O}_X}
\mcal{P}_X \otimes_{\mcal{O}_X} \mcal{M}_j)|_{U_{\bsym{i}}} \Bigr) 
\\
& \qquad \qquad \to \
\Omega(\bsym{\Delta}^l_{\K}) \hatotimes{}
(\Omega_X \otimes_{\mcal{O}_X}
\mcal{P}_X \otimes_{\mcal{O}_X} \mcal{N})|_{U_{\bsym{i}}} \big) 
\end{aligned} \] 
indexed by $\bsym{i}$ and satisfying the simplicial relations. 
Now use Lemma \ref{lem3.2} to obtain $\bsym{\sigma}^*(\phi)$.
Properties (i-iii) follow from Lemma \ref{lem6.6}.
\end{proof}


\end{document}